\documentclass[10pt,a4paper]{article}
\usepackage[utf8]{inputenc}\usepackage{pdfsync}\usepackage{float}
\usepackage[T1]{fontenc}
\usepackage{amsmath}
\usepackage{amsfonts}
\usepackage{amssymb}
\usepackage{amsthm}
\usepackage{xcolor}
\usepackage{bbm,graphicx}
\usepackage[left=2cm,right=2cm,top=2cm,bottom=2.5cm]{geometry}
\usepackage{authblk}
\usepackage{subcaption}
\usepackage{hyperref}
\newtheorem{thm}{THEOREM}%[section]

\newtheorem{cor}[thm]{COROLLARY}
\newtheorem{lem}[thm]{Lemma}
\newtheorem{definition}{DEFINITION}
\newtheorem{example}{EXAMPLE}
\newtheorem{proposition}[thm]{PROPOSITION}
\newtheorem{remark}{REMARK}
\def\mA{{\mathcal A}} \def\mB{{\mathbb B}} \def\mC{{\mathcal C}}
\def\beXa{\begin{example}} \def\eeXa{\end{example}}
\def\eeD{\end{definition}} \def\beD{\begin{definition}}
\def\beR{\begin{remark}} \def\eeR{\end{remark}}
\def\beL{\begin{lem}} \def\eeL{\end{lem}}
\def\beP{\begin{proposition}} \def\eeP{\end{proposition}}
\def\beC{\begin{cor}} \def\beT{\begin{thm}}
  \def\eeT{\end{thm}}
\def\eeC{\end{cor}}
 \def\lc{loading coefficient}\def\mbe{may be expressed in terms of } \def\mbr{may be reduced to }
\def\itm{it may be checked that }
\providecommand{\pr}[1]{\left(#1\right)} %(.)
\providecommand{\pp}[1]{\left[#1\right]} %[.]
\providecommand{\set}[1]{\left\lbrace#1\right\rbrace} %{.}
 \def\ssec{\subsection}

    \def\Tq{\tilde{q}} \def\Tl{{\tilde{\lambda}}}     \def\fun{function } \def\funs{functions } 
\def\GS{Gerber-Shiu } \newcommand{\bff}[1]{{\mbox{\boldmath$#1$}}}
 \def\wk{well-known } 
   \def\para{parameter} \def\paras{parameters} 
  \def\rps{ruin probabilities}  
 \def\wk{well-known} 

\def\G{\Gamma}
  
 \def\procs{processes} \def\proc{process}
\def\Rui{\Psi} 
 \def\fp{first passage }
\def\ovl{\overline}
\def\lev{L\'evy }   \def\rp{ruin probability} \def\Oth{On the other hand, }
\def\dd{draw-down } \def\sn{spectrally negative }
\long\def\symbolfootnote[#1]#2{
\begingroup
\def\thefootnote{\fnsymbol{footnote}}\footnote[#1]{#2}
\endgroup}
\def\fn{\symbolfootnote}

\def\I{\infty} \def\Eq{\Leftrightarrow}\def\vb{\vec \beta}\def\vo{\bff 1}
  \def\T{\widetilde}
\def\CL{Cram\`er-Lundberg }  
 \def\PH{phase-type }  
\def\BEN{\begin{enumerate}}  \def\BI{\begin{itemize}}
\def\EEN{\end{enumerate}}   \def\EI{\end{itemize}} \def\im{\item} \def\Lra{\Longrightarrow}  \def\eqr{\eqref}  
\def\no{\nonumber} %\newcommand{\be}{\begin{equation}}
\def\mG{{\mathcal G}}
\def\g{\gamma}  \def\d{\delta}   
  \def\th{\theta} 
\def\e{\epsilon} 
\def\k{\kappa} \def\l{\lambda}  %\def\W{W}
\def\Lm{L\'evy measure }   \def\resp{respectively}  \def\SLG{Shreve, Lehoczky and Gaver}  \def\PK{Pollaczek-Khinchine } \def\LZ{Lokka-Zervos alternative}
\def\lm{\nu}  \def\ith{it holds that } \def\sf{scale function}
   \def\s{\sigma}  
  \def\Fq{\Phi_{q}} 
  \def\satd{satisfied } \def\Z1{Z_{\q,1}}
  \def\bc{\begin{cases}
  } \def\gene{generalization }    \def\cP{compound Poison model }
\def\ec{\end{cases}} \def\expo{exponential } 
   \def\for{\forall}
  \def\beq{\begin{eqnarray}} \def\eeq{\end{eqnarray}}
   \def\be{\begin{equation}} \def\ee{\end{equation}}
\def\bea{\begin{eqnarray*}} \def\rei{reinsurance } \def\mbw{may be written as }
  \def\prf{{\bf Proof:} } \def\satg{satisfying } 
\def\eea{\end{eqnarray*}} \def\la{\label}
 \def\LT{Laplace transform} \def\BM{Brownian motion}
 \def\q{q} \def\R{{\mathbb R}}
\def\le{\left} \def\ri{\right}
\def\nne{nonnegative } \def\fr{\frac}    \def\ta{T_{-a}}  \def\tb{T_{b,+}}  
\def\tz{T_{0,-}} \def\sec{\section} \def \fe{for example } \def\wlo{w.l.o.g. }  \def\ith{it holds that }
 
 \def\deF{de Finetti } \def\app{approximation}

 \def\cP{compound Poisson }

  \def\expo{exponential } 
\def\prob{problem}   
 \def\How{However, } \def\H{\widehat }
\def\mC{\mathcal C}   
\def\1{\mathop{\rm 1\!\!I}\nolimits}

\def\prob{problem}   \def\str{strong Markov property}

  \def\Itf{It follows that }    
   
\def\me{matrix exponential jumps}
 \def\sp{spectrally positive } 
 \newcommand{\md}{\mathrm d}   \def\divs{dividends }\def\unk{unknown}\def\mbs{may be shown }
\newcommand{\red}{\textcolor[rgb]{1.00,0.00,0.00}}

   \def\Vb{V^{b]}} \def\Vkb{V_k^{b}}  
     
   \def\Fr{Furthermore, }  \def\app{approximation}
       \def\mA{\mathcal A} \def\mB{\mathcal B} \def\mC{\mathcal C}   \def\ubd{unbounded variation case} 
\def\bd{bounded variation case} \def\div{dividend}  \def\Pd{Pad\'e approximation } \def\Pds{Pad\'e approximations}\def\tPd{two point Pad\'e approximation}  
    \def\Lm{L\'evy measure } \def\wkt{well-known and easy to check that } \def\expoc{exponential claims}
    \def\sats{satisfies}   \def\fun{function }  \def\funs{functions} \def\thr{therefore} \def\wf{we find that }  \def\resp{respectively} \def\proc{process}     \def\ci{capital injections} \def\C{I} \def\D{D}  
 \def\part{\partial }    \def\bep{\begin{pmatrix}} \def\eep{\end{pmatrix}}
       \def\wrt{with respect to } 
           \def\c{c}  \def\eddc{expected discounted dividends minus capital injections} \def\per{perturbed } \def\Ito{It turns out that }  \def\saty{satisfy } \def\exc{excursion theory} \def\pt{provided that }
   
       \def\CLp{\CL process  }  \def\fno{from now on}   \def\expoj{exponential jumps}  \def\Oth{On the other hand, }\def\snL{spectrally negative L\'evy processes} \def\GSf{\GS function}\def\deV{de Vylder approximation}\def\wmh{we must have }\def\vRui{\vec \Rui}\def\frq{\frac q{\Fq}}
   
   \def\capi{capital injections}   \def\Nt{Note that }
\providecommand{\pr}[1]{\left(#1\right)} %(.)
\providecommand{\pp}[1]{\left[#1\right]} %[.]
\providecommand{\set}[1]{\left\lbrace#1\right\rbrace} %{.}
 \def\si{\fr{\s^2}2}\def\re{\textcolor{red}}\def\mA{\mathcal A} \def\mB{\mathcal B} \def\mC{\mathcal C} \def\mD{\mathcal D}\def\sm{smooth fit}\def\strf{straightforward}
 
 \def\cm{completely monotone}
\newcommand{\figu}[3]{
\begin{figure}[H]
\begin{center}
{\includegraphics[width=13 cm, height=8 cm]{#1}}
\end{center}
\vspace{-0.2cm}
\caption{\hspace{0.25cm}#2\label{f:#1}}
\end{figure}
}

\begin{document}
\title{Optimizing dividends and capital injections limited by bankruptcy, and  practical approximations for the  Cram\'er-Lundberg process
}
%\title{Are the  ``de Vylder-type" and Pad\'e  approximations   efficient for optimizing dividends for the \CL process?}
%\author{F. Avram, D. Goreac, R. Adenane, U. Solon }

 \author[1]{F. Avram \footnote{Corresponding author. E-mail:
     Florin.Avram@univ-Pau.fr}}
 \author[2,3]{D. Goreac}
 \author[4]{R. Adenane}
 \author[1,5]{U. Solon}
 \affil[1]{Laboratoire de Math\'{e}matiques Appliqu\'{e}es, Universit\'{e} de Pau, France 64000 }
 \affil[2]{LAMA, Univ Gustave Eiffel, UPEM, Univ Paris Est Creteil, CNRS, F-77447 Marne-la-Valleée, France}
 \affil[3]{School of Mathematics and Statistics, Shandong University, Weihai, Weihai 264209, PR China}
 \affil[4]{Département des Mathématiques, Université Ibn-Tofail, Kenitra, Maroc 14000}
 \affil[5]{Institute of Mathematics, University of the Philippines Diliman, Philippines}

\maketitle

\begin{abstract}
The recent papers Gajek-Kucinsky(2017) and Avram-Goreac-Li-Wu(2020) investigated  the  control problem  of  optimizing dividends when limiting capital injections stopped upon bankruptcy. % is taken into consideration.
The first paper works under the  spectrally negative L\'evy model; the second works under the  Cram\'er-Lundberg model  with exponential jumps,  where the results are considerably more explicit.  The current paper has three purposes. First, it illustrates the fact that quite reasonable approximations of the general problem may be obtained using the particular exponential case studied in Avram-Goreac-Li-Wu(2020).  Secondly, it extends the  results   to the case when a final penalty $P$  is taken into consideration as well besides a proportional cost $k>1$ for \ci. This requires amending the ``scale and Gerber-Shiu functions" already introduced in Gajek-Kucinsky(2017). Thirdly, in the exponential case, the
results will be made even more explicit by employing the Lambert-W function. This tool has particular importance in computational aspects and can be employed in theoretical aspects such as asymptotics.
\end{abstract}
\tableofcontents

\section{Introduction %Pad\'e and two-point Pad\'e approximations, with low order examples
\label{s:low}}

This paper concerns the approximate optimization of a new type of boundary  mechanism, which emerged recently in the actuarial literature \cite{APY,Gaj,AGLW},
in the context of  the optimal control of dividends and capital injections.

\indent {\bf The model}. Consider the \sn \lev risk model:
\begin{equation*}
\label{CLp}
X_t=x+c t-\bar{\xi}_t,
c \geq 0,
\end{equation*}
where $\bar{\xi}_t$ is a \sp \lev \proc, with \Lm $\lm(x)dx$.
The classic example is that of the \per \CL risk model with $$\bar{\xi}_t=\sum_{i=1}^{N_t}\xi_i + \s B_t,$$ where $B_t$ is a
\BM, where $N_t$ is an independent Poisson process of intensity $\lambda>0$, and  $\pr{\xi_i}_{i\geq 1}$ is a family of i.i.d.r.v.  whose distribution, density and   moments are denoted respectively by $F,f,m_i, i \in \{1, ... \}$. % and $B_t$ is an independent \sta \BM.
\iffalse
The space is then endowed with the natural right-continuous, completed filtration
 $\mathbb{F}$  satisfying the usual assumptions of right-continuity and completeness.
\fi
\Fr
\begin{itemize}
\item the process is modified by {\bf dividends and capital injection}:
$$\pi:=\pr{\D,I}\Rightarrow X_t^\pi:=X_t-\D_t+\C_t,$$
where $\D,I$ are adapted, non-decreasing and càdlàg processes, with   $\D_{0-}=\C_{0-}=0$;
\item the  first time when we do not bail-out to positive reserves $\s^{\pi}_{0}:=\inf\set{t>0:\ X_{t-}^\pi-\bigtriangleup \bar{\xi}_t+\bigtriangleup \C_t<0}$ is called {\bf bankruptcy/absolute ruin};

    \item prior to  bankruptcy, dividends are limited by the available reserves: $\bigtriangleup \D_t:=\D_{t}-\D_{t-}\leq X_{t-}^{\pi}-\bigtriangleup \bar{\xi}_t+\bigtriangleup \C_t, \textnormal{ where }\bar{\xi}_t:=\sum_{i=1}^{N_t}\xi_i$. The set of ``admissible" policies \satg this constraint will be denoted by
        $\tilde{\Pi}(x)$.
\end{itemize}

The objective is to maximize the {\bf profit}:
$$J_x^\pi:=\mathbb{E}_x\pp{\int_{\pp{0,\s^{\pi}_{0}}}e^{-qt}
\pr{d\D_s-kd\C_s}-P e^{-q \s^{\pi}_{0}}}, k\geq 1, P \geq 0.$$
The {\bf value function} is  $$V\pr{x}:=\sup_{\pi\in\tilde{\Pi}(x)}J_x^\pi,\ x\in \mathbb{R}.$$

{\bf Motivation.} The recent papers \cite{Gaj,AGLW}
 investigated  the  above control problem  of  optimizing dividends
 and  capital injections for \procs\ with jumps, when bankruptcy is allowed as well.
 The second paper works under the  Cram\'er-Lundberg model  with \expoj, while the first  works under the  spectrally negative L\'evy model, allowing also for the presence of Brownian motion and infinite activity jumps.
 \Ito
 the optimal policy belongs to the class of $(-a,0,b), a>0, b \geq 0$ ``bounded buffer policies", which
consist in allowing only capital injections smaller than a given $a$ and declaring bankruptcy at the first time when the size of the overshoot below 0 exceeds $a,$ and  in paying dividends when the reserve reaches an upper barrier $b$. These will briefly be described as $(-a,0,b)$ policies \fno.
\Fr the optimal $(a^*,b^*)$ are the
        roots of  one variable equations with explicit solutions related to the Lambert-W(right) function (ProductLog in Mathematica).

\iffalse
 In this paper, we revisit
the computation of the value of  $(-a,0,b)$ policies; this requires  conditioning at the first exit from $[0,b]$, and also at the draw-down time of the first capital injection, starting from $b$. At the same time,  we allow    for an extra final penalty $P$, which does not pose problems at the level of computing the value of  $(-a,0,b)$ policies. Note though that optimality is not  investigated below.

In a different direction, we exploit the fact that   results are considerably more explicit in \cite{AGLW}, and even more below, where we express them using the LambertW function,
to offer exponential approximations for the general case.
\fi

Below, our goal is to  show numerically that  exponential approximations provide quite reasonable results  (as the \deV\ provides for the ruin problem).  We will focus  in our examples  on the case of \me (which are known to be dense in the class of general \nne jumps, with even error bounds for \cm\ jumps being available \cite{vatamidou2014accuracy}), for two reasons. One is in order to highlight  certain exact equations which are similar to their exponential  versions,  and
 which may  at their turn be used to produce even more accurate approximations in the future, and, secondly,   since  numerical Laplace inversion for this class may easily  tuned to have arbitrarily small errors.
 % which are both general and practical for computations.

 {\bf History of the problem}:
The case of no capital injections (also characterized by $k=\infty$ or absorption below 0) is the dividend problem posed by De Finetti \cite{DeF,gerber1969entscheidungskriterien} where dividends are paid above barrier $b^*$ and $a^*=0$ is imposed. ``The challenge is to find the right compromise
between paying early in view of the discounting or paying late in order not to reach ruin too early and thus profit from the  positive safety loading for a longer time" \cite{albrecher2020optimal}.

Forced injections and no bankruptcy at $0$ (also characterized by a reflection at $0$) is studied in Shreve \cite{SLG} where dividends are paid above barrier $b^*$ and $a^*=\infty$ is imposed.

From Lokka, Zervos, \cite{LZ} we know that in the  \BM\ case, it is optimal to either always inject, if  $k\leq k_c$, for some critical cost $k_c$ (i.e. use Shreve), or, stop at $0$ (use De Finetti). We propose to call this the \textbf{Lokka, Zervos alternative}. The ``proof" of  this alternative starts  by largely assuming it via a {heuristically justified border Ansatz} [LZ08, (5.2)]:
$\max\set{-V(0),\ V'(0)-k}=0 \Lra \text{ either } V(0)=0 \text{ or } V'(0)=k.$

Extensive literature {\bf on SLG  forced bailouts (no bankruptcy)} can be found at Avram et al., (2007) \cite{APP}, Kulenko and Schmidli, (2008) \cite{kulenko2008optimal}, Eisenberg and Schmidli, (2011) \cite{eisenberg2011minimising}, P\'erez et al., (2018) \cite{perez2018optimal}, Lindensjo,  Lindskog (2019) \cite{lindensjo2019optimal}, Noba et al., 2020 \cite{noba2020bailout}.

Articles \cite{Gaj,AGLW} are the only papers which relate declaring bankruptcy to the  size of  jumps, with general and exponential jumps, \resp. \cite{Gaj} deals  also with the presence of Brownian motion and infinite activity jumps, by conditioning at the first draw-down time; the optimality proof  is %(naturally)
 quite involved.

In \cite{AGLW}, it is also shown that neither $V(0)=0$ \text{ nor } $V'(0)=k$ are possible: the \LZ\ disappears, but another interesting alternative holds.
Above   a certain critical $k_c$  the optimal \divs\ barrier switches from strictly positive to $0$, and  $k_c$
is related in \eqr{kc} to the Lambert-W \fun.

 The results of \cite{Gaj,AGLW} may be divided in three parts:
 \BEN \im Compute the value of bounded buffer policies. The key result is \eqr{struct} below,
 an explicit determination of the objective $J_0=J_0^{a,b}$, which allows optimizing it numerically.

\beR Computing the value function   is considerably simplified  by the use of  \fp  recipes   available for spectrally negative L\'evy processes \cite{AKP,Kyp,KKR,AGV},  which are built around   two ingredients: the  $W_q $ and $Z_q$ $q$-scale functions, defined respectively  for $x \geq 0, q \geq 0$ as:\BEN \im the inverse Laplace transform  of $\fr 1{\k(s)-q}$, where $\k(s)$ is the Laplace exponent (which characterizes a \lev \proc) and \im
$Z_q(x)=1+q\int_0^xW_q(y)dy$\EEN
 -- see the papers \cite{Suprun,Ber,AKP} for the first appearance of these functions.  The name $q$-scale/harmonic
  functions  is justified by the fact that these functions are harmonic for the
process $X$ killed upon entering $(-\I, 0)$, in the sense that
$$\{e^{- q \min[t,T_0]}\;  W_q(X_{\min[t,T_0]}),   e^{- q \min[t,T_0]}\; Z_q(X_{\min[t,T_0]})\}, t \geq 0$$
are martingales, as shown in \cite[Prop. 3]{Pisexit} (in the case of $Z_q$, there is also a {\bf penalty} of $1$ at ruin, generalizing to other  penalties produces the so-called \GSf).

\eeR

\im Equations determining candidates for the optimal
$a^*,b^*$ are  obtained
 by differentiating  the objective (which is expressed in terms of the scale functions
$W_q,Z_q$), and  the   optimal pair
    $(a^*,b^*)$ is identified. As a result, the critical  $k_c$ is related in \eqr{kc} to the  Lambert-W function.
    \im The optimality of the $(-a^*,0,b^*)$ policy is established.
 \EEN

Note that the last  step is   quite non-trivial and is achieved by different methods in \cite{Gaj} and \cite{AGLW}. The latter paper starts by formulating a (new)  HJB equation associated to this stochastic control problem -- see \eqr{HJB}.

\beR  The objective may be optimized numerically using the first step only (the equation \eqr{J0} for $J_0^{a,b}$).
\eeR

 Exponential approximations may also be used, which  are  similar in spirit with the de Vylder-type approximations. Recall that the philosophy of the \deV\  is to approximate a \CLp by a simpler \proc\ with \expoj, with cleverly chosen exponential rate $\mu$,  and the parameters $\l,c$ may also be modified, if one desires to make the approximations exact at $x=0$ --see \fe \cite{AHPS} for  more details)

 The efficiency of the de Vylder \app\ for approximating \rps\ is well documented \cite{de1978practical}. The natural question of whether this type of techniques may work for other objectives, like \fe\ for optimizing dividends and/or \rei\ was already discussed  in \cite{hojgaard2002optimal,dickson2005optimal,beveridge2007optimal,GSS,AHPS}. In this paper, following on previous works \cite{avram2011moments,AP14,ABH}, we draw first the attention to the fact that  we have not one, but three  de Vylder-type approximations for $W_q(x)$ (as for the \rp). The best \app\ in our experiments when the \lc\ $\th$ is large turn out to be the classic \deV. However, for approximating near the origin, the two point \Pd  which fixes both the values $W_q(0)=\fr 1 c, W_q'(0)=\fr {q + \l}{c^2}$ works better. The end result here is simply replacing the inverse rate $\mu^{-1}$     by $m_1, $ in the formula for the \sf\ of the \CLp\ with \expoj. In between $x=0$ and $x \to \I$,  the winner is sometimes  the ``Renyi approximation" which replaces the inverse rate  by $\fr {m_2}{2 m_1}$, and modifies $\l$ as well (for the \deV, $\mu^{-1}$ is replaced  by $\fr {m_3}{3 m_2}$, and both $\l,c$  are modified).

 We end this introduction by  highlighting in  figure \ref{f:ZZ} the fact  that for \expoj, the limited capital injections objective function $J_0$  given by \eqr{Estima*} for arbitrary $b$ but optimal $a=s(b)$ (via a complicated formula)  improves
the value function \wrt de Finetti and \SLG, for any $b$.

\figu{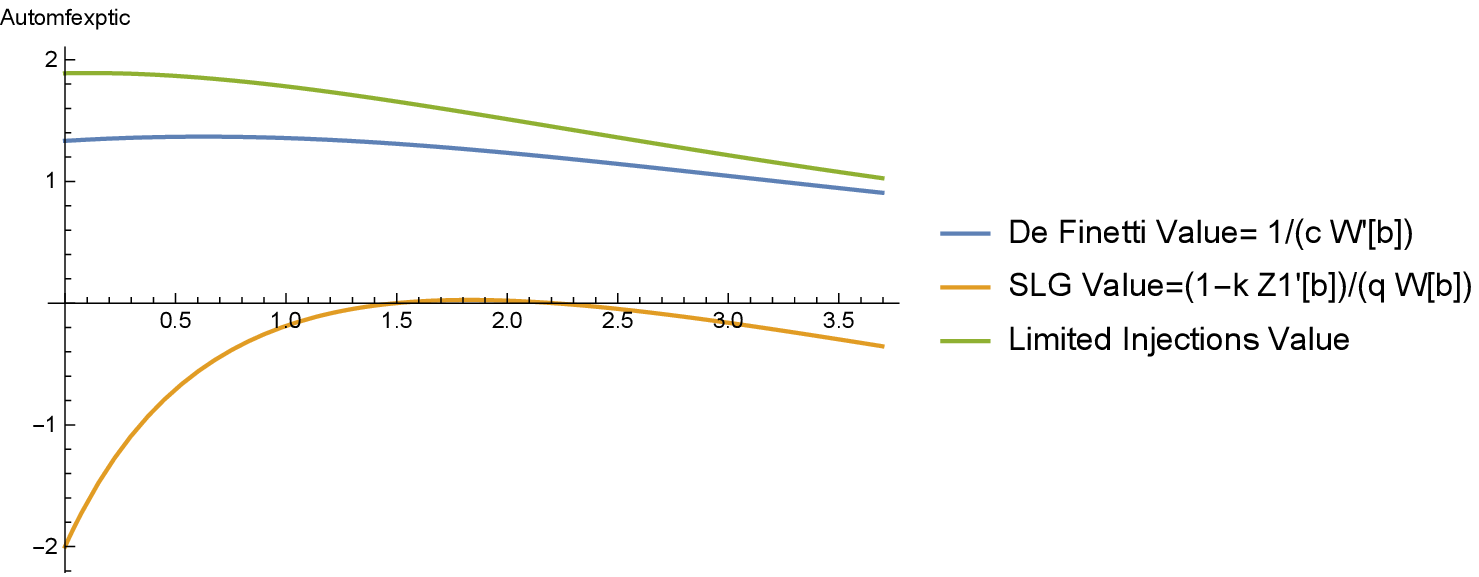}{The value function $J_0$  given by \eqr{Estima*}, for arbitrary $b$ but optimal $a$. The inequality observed is a  consequence of the properties of the Lambert function. The improvement \wrt de Finetti is considerable, of $0.382292 \%$ (the SLG approach is not competitive in this case).  Note also that the optimal barrier $b=0.109023$ is smaller than the de Finetti and SLG optima of $0.626672, 1.82726$ respectively.}{.7}

{\bf Contents and contributions}.
Section \ref{s:AG} offers a conjectured profit formula for $(-a,0,b)$ policies, where we include also a final penalty $P$. The theoretical result of the section \ref{p:cost} revisits \cite{Gaj,AGLW} by linking the two formulations together and emphasizing the impact of the bankruptcy penalty $P$ (via the scale function $G$). Its proof is beyond the scope of the present (already lengthy enough) paper and it can be inferred from either one of \cite{Gaj} and \cite{AGLW} through a three step argument:
\begin{enumerate}
\item express the cost by conditioning on the reserve ($J_x$) starting from $0 \leq x \leq b$ hitting either 0 or b;
\item get a further relationship on costs $J_b$ and $J_0$ by conditioning on the first claim;
\item finally, mix these conditions together in order to obtain the explicit formula for $J_x$.
\end{enumerate}

We also wish to point out (in section \ref{s:HJBsys}) the link to an appropriate HJB variational inequality equation\ref{HJB} and the definition \ref{d:HJB} specifying the action regions and their computation starting from
the regimes in the HJB system. To the best of our knowledge, this is new and the preliminary studies conducted on more complicated problems (involving reinsurance and reserve-dependent premium) seem to reinforce the relevance of this tool.

In  section \ref{s:me} we  provide an alternative matrix exponential form of the exact cost, in  the case of \me.

An explicit determination of $a^*,b^*$ and an  equity cost  dichotomy  when dealing with exponential jumps are given in section \ref{s:eqc}, taking also advantage of properties of the Lambert-W function, which were not exploited before. The two main novelties of the section are:
\begin{itemize}
\item emphasizing the computations of the optimal buffer/barrier (from \cite{AGLW}) in relation with the scale-like quantities appearing in \cite{Gaj};
\item making explicit use of the (computation-ready) Lambert-W function to describe the dependency of optimal $a^* b$ (in equation \ref{Estima*}) and of the dichotomy-triggering cost $k_c$ in equation \ref{kc}.
\end{itemize}
Again, a further novelty is the presence of the bankruptcy cost $P$.

Section \ref{s:DV} reviews, for completeness,  the \deV-type \app s.  Section \ref{s:DVr} recalls, for warm-up, some of the oldest exponential approximations for \rps. Section \ref{s:DVW} recalls in Proposition \ref{p:AH}, following \cite{AP14,AHPS}  three approximations of the scale function  $W_q(x)$\fn[4]{essentially, this is the  ``dividend function with fixed barrier", which had been also extensively studied in previous literature before the introduction of $W_q(x)$}, obtained by approximating its \LT. These amount finally to replacing our process by one with exponential jumps \ and cleverly crafted parameters based on the first three moments of the claims.

In section \ref{s:ci}, we consider particular examples and obtain very good approximations for two fundamental objects of interest: the growth exponent $\Fq$ of the scale function $W_q(x)$, and the (last) global minimum of $W_q'(x)$, which is fundamental in the de Finetti barrier problem. Proceeding afterwards to the problem of dividends and limited capital injections, concepts in section \ref{s:eqc} are used to compute a straightforward exponential approximation based on an exponential approximation of the claim density, and a new ``correct ingredients approximation" which consists of plugging into the  objective function \eqr{J0} for exponential claims the exact "non-exponential ingredients" (scale functions and, survival and mean functions) of the non-exponential densities. Both methods are observed to yield reasonable values in approximating the objective.

This leads us to our conclusion that from a practical point of view, exponential approximations are typically sufficient in the problems discussed in this paper.

%Section \ref{s:Ma} gives some idea of the programs we used, which are available upon request

\iffalse
Finally, section \ref{dev} recalls the derivation of some  de Vylder-type approximations, including the  original derivation of the Renyi and de Vylder approximations using  process cumulants, in section \ref{derdev}.

To prevent this paper form becoming too long, we decided to postpone for the future the investigation of the performance of the approximation \eqr{J0} on  some non-matrix exponential favorites of the  statistical modeling like   gamma (including $\chi^2$), Pareto, Weibull, Mittag-Leffler, and beta.
In that case, our approximation must be tested against the exact formula of \cite{Gaj}.
\fi

% section \ref{s:two} recalls reviews some aspects of the  asymptotic behavior of the \sf, which are crucial in evaluating numerically the \perf\of our approximations.

\section{The cost function of $(-a,0,b)$ policies, for the \sn \lev case \la{s:AG} %\la{s:cost}
}

In this section, we allow  $\bar{\xi}_t$ to be a \sp \lev \proc, with a \Lm admitting a density $\lm(dy)=\nu'(y) dy$.
The simplest example is that of the \per \CL risk model with $$\bar{\xi}_t=\sum_{i=1}^{N_t}\xi_i + \s B_t,$$ where $N_t$ is a Poisson process of intensity $\lambda>0$, $\pr{\xi_i}_{i\geq 1}$ is an independent family of i.i.d.r.v. with density $f(y)$, and   $B_t$ is an independent
\BM.

 We revisit here the   problem of optimizing  the value of "bounded buffer $(-a,0,b)$ policies", following \cite{Gaj,AGLW} (in order  to relate the results, one needs to replace  $\g$ in the  objective of \cite{Gaj} by  $1/k$), while taking into account also the bankruptcy penalty $P$.

 An  important role in the results  will be played by the {\bf expected scale after a jump}
\be C(x)=\int_0^x  W_q(x-y)\ovl \lm (y) dy=c {W_q(x)} -Z_q(x)+\si W_q'(x),
\la{C}\ee
where $\ovl \lm (y)=\int_y^\I \lm(u) du$ is the tail of the \lev measure and $\s$ is the Brownian volatility (the identity above follows easily from the $q$-harmonicity of $Z_q$, after an integration by parts of the convolution term and a division by $q$).

The problem of limited reflection requires introducing a  new "scale function $S_a(x)$ and Gerber-Shiu function $G_{a,\s}(x)$"--
see Remark \ref{r:wr} for further comments on this terminology:
\be \la{RG} \bc S_a(x)=Z_q(x)+  C_a(x), C_a(x) = \int_0^x W_q(x-y) \;   \ovl \lm(a+y) \; dy\\
G_{a,\s}(x) =
G_a(x) + k \si W_q(x) \ec \ee
where
\bea && G_a(x)=  \int_0^x  W_q(x-y) \pr{k \;  m_a(y)+ P \ovl \lm(a+y) }   dy :=k M_a(x)+ P C_a(x), \\ &&m_a(y)=\int_0^a z \nu(y+z)dz.  \eea

\beXa
With
\expoj\ and possibly $\s>0$, using the  identities $$\ovl \lm(y+a)=e^{-\mu a} \ovl \lm(y), m_a(y)= \l e^{-\mu y} m(a), m(a)=\int_0^a y \; \mu  e^{-\mu y} d y=\frac{1-e^{-\mu a }}{\mu }-a e^{-\mu a },$$ we find that the  functions  \eqr{RG}   are expressible {\bf as products} of $C(x)$ and the {\bf survival} or {\bf mean} function of the jumps:
\be \bc  C_a(x)=
C(x)e^{- \mu  a}=C(x) \ovl F(a), \; \ovl F(a)=1-F(a)\\S_a(x) = Z_q(x) +e^{-\mu a} C(x)\\G_a(x) = \pr{k m(a) +P \ovl F(a)} C(x)
\ec \ee
(\cite{Gaj} use $s_c,r_c$, instead of $M_a(x):=\int_0^x  W_q(x-y) \;  m_a(y) \; \; dy,C_a(x)$, \resp). When $P=0=\s$,
these reduce to quantities in \cite{AGLW}.

The   formulas above  will be used below as a heuristic approximation in non-\expo\ cases.

\eeXa

\beR Note that
\be \la{ini} C_a(0)=0, G_a(0)=0, S_a(0)=1, C(0)=0,  C'(0)=\bc \fr \l c &\s=0\\0 &\s>0\ec, %C''(0)=\bc \fr \l {c^2}\pr{\l +q - \mu c} &\s=0\\-\frac 2 {\sigma^2}(\frac {2 c^2} {\sigma^2}+q) &\s>0\ec
\ee
 and that $C(x),G_a(x),S_a(x)$ are increasing functions in $x$.
\eeR
 We state now a %conjectured 
 \gene\ of \cite[Thm. 4]{Gaj} for the value function  $J_0^{a,b}$ of $(-a,0,b)$ policies, in terms of  $S_a(x), G_a(x)$.    In the \CL case illustrated below, the proof is \strf,  following \cite{AGLW}. In the other case, one needs to adapt the proof of \cite{Gaj}. %Note that the explicit  formula \eqr{struct} may be optimized numerically.

\beT  {\bf Cost function  for $(a,b)$ policies} \la{p:cost}
 For a  \snL, let
$$J_x=J^{a,b}(x):=
\mathbb{E}_x\pp{\int_0^{\ta}e^{-qt}\pr{\md \D_t -k \;\md\C_t}-P e^{- q \ta}}  $$
denote the \eddc\ associated to policies consisting in paying capital injections  with proportional cost $k\geq 1$, provided that the severity of ruin is smaller than $a>0$, and paying dividends as soon as the  process reaches some upper level $b$.
Put \bea G_{a,\s}(x)=  G_a(x) + k \si W_q(x).\eea

 Then, \ith
 \be J_x= \bc  G_{a,\s}(x) +J_0^{a,b} S_a(x)= G_{a,\s}(x) +\fr{1- G_{a,\s}'(b)}{S_a'(b)} S_a(x),  &x \in[0, b]\\k x+J_0^{a,b}
&x \in[-a, 0]\\ 0 & x \leq -a \ec. \la{struct}\ee

\eeT

\beR \la{r:wr}

The first equality in \eqr{struct} will be easily obtained by applying  the \str\ at the stopping time $T=\min[\tz,\tb]$, but it still contains the \unk\ $J_0$.

This relation suggests a definition of the scale $S_a$ and the \GSf\ $G_{a,\s}$,  as the coefficient   of $J_0$ and  the part    independent of $J_0$, \resp.

This equality  is also equivalent to
\be \la{wr}
J_0^{a,b}= \fr{J_x- G_{a,\s}(x)}{S_a(x)}=\fr{1- G_{a,\s}'(b)}{S_a'(b)},\ee
which suggests
another {\bf analytic definition of the scale and \GSf} corresponding to an objective $J_x$ which involves reflection at $b$.

 The \funs\ $S_a(x), G_a(x)$ \mbs to stay the same for problems
which require only modifying the boundary condition at $b$, like the problem of \ci\ for the process reflected at $b$, or the  problem of \divs\ for the process reflected at $b$, with
proportional retention $k_D$ (this is in coherence with previously studied problems).

\eeR

\beC
Let us consider the \CL\ setting without diffusion (i.e. $\s = 0$),
 For fixed  $k \geq 1, b\geq 0$, the optimality equation $\fr{\partial}{\partial a} J_0^{a,b}=0$  \mbw \be \la{paag}
J_0^{a,b}=k a -P  \Eq J_{-a}^{a,b}= -P. \ee

\eeC
\beR The first equality in \eqr{paag} provides a relation between the objective $J_0$ and the variable $a$;  the second recognizes this as the smooth fit equation $J_{-a}=0$.
\eeR
\prf
Recalling the expressions of $J_0^{a,b}$, $G_a(x)$,  in  \eqr{wr},  in \eqr{RG},
and from \cite[Lem. A.4]{Gaj} $$M_{a'}(x)= -a  C_{a'}(x),$$ where
 $C_{a'}(x),M_{a'}(x)$ denote derivatives with respect to the subscript $a$,
 Whenever $b>0$, if $a$ achieves the maximum in $J^{a,b}_0$, it is straightforward (think of the economic interpretation) that $a$ achieves the maximum of $a\mapsto J^{a,b}_x$ for every $x\in [0,b]$. Therefore, 
we find
\bea
\begin{aligned}
& \fr{\partial}{\partial a} J_0^{a,b}=0 \Eq  J_0^{a,b}=\frac{- G_{a'}(x)}{C_{a'}(x)}= \frac{-k M_{a'}(x) - P C_{a'}(x)}{C_{a'}(x)}= k a -P\\
& \Eq  J_{-a}^{a,b}=J_0^{a,b}- ka \Eq J_{-a}^{a,b}= -P.
 \end{aligned}
\eea

\qed
\iffalse
\beR Note that
\be (S(x),G_{k,P}(x))=\bc (W_q(x),0)&s=0\\(Z_q(x),\fr{C(x)}\mu )\re{=?}(Z_q(x),\Z1(x)- \si W_q(x))&s=\I\ec, \ee
where the last identity
\be \fr{C(x)}\mu +  \si W_q(x)=\Z1(x) \la{CZ1}\ee
is specific to this case
and generalizes an identity  established in \cite{AGR}.

Our formula \red{?interpolates } between the de Finetti and \SLG cases:
\be  \la{obp} J^{s,b}(x)=\bc \Vb(x)=\fr{W_q(x)}{W_q'(b)}&s=0\\\Vkb(x)=k   {\Z1(x)} + Z_q(x) \fr{1-k \Z1'(b)}{Z_\q'(b)}&s=\I\ec.\ee

\eeR

\fi

%\input{Pr1}

\subsection{The HJB System} \label{s:HJBsys}

The  optimality proof in
 \cite{AGLW} is based on showing  that the function $J_x$ \eqr{struct} with $a^*,b^*$ defined in \eqr{J0a},
 \eqr{strEqb} is the {minimal AC-supersolution of the   %, everywhere except at $0$ and at $-a^*$,
  HJB system
   \begin{equation}
\label{HJB}
\left\lbrace
\begin{split}
&\max\set{ {H}\pr{x, {V}, {V}'(x)},1- {V}'(x), {V}'(x)-k}=0,\ \forall x\in\mathbb{R}_+\\
&\max\set{ {V}'(x)-k,-P-V(x)}=0,\ \forall x\in\mathbb{R}_-
\end{split}\right.,
\end{equation}
where the Hamiltonian $ {H}$ is given by
\begin{equation}
\label{Hamiltoniantilde}
  {H}\pr{x,\phi,v}:=c v+\lambda\int_{\mathbb{R}_+}\phi(x-y)\mu e^{-\mu y}-\pr{q+\lambda}\phi(x).
\end{equation}

To discuss  \eqr{HJB}, it is useful to introduce the concept of {\bf dividend-limited injections
 strategies} and {\bf barrier strategies}. The following are also valid for its generalizations to mixed singular/continuous controls taking into account reinsurance:
\beD \label{d:HJB} Dividend-limited injections strategies
are  stationary strategies where the dividends are paid  according to a partition of the state space $\R$ in five sets  $ \mA, \mB,\mC, \mC_0,\mD$ as follows:
\BEN \im  If the surplus is
in $\mA$ (absolute ruin), bankruptcy is declared and a penalty $P$ is paid;\im  If the surplus is
in $\mB$ bailouts/\capi\ are used for bringing the surplus to the closest point of $\mC$.
\im  If the surplus is
in the open set $\mC$ (continuation/no action set), %which is an  interval where
no controls are used. \im
 If the current surplus is in
$\mC_0\subset \mathcal{D}$ (these are upper-accumulation points of $\mathcal{C}$),
dividends are paid at a positive rate,  in order to keep the surplus process from moving.
\im
 If the current surplus is in
$\mD$, a positive amount of money is paid as
dividends in order to bring the surplus process to $\mC_0$.

 {\bf Barrier strategies} are stationary strategies for which $\mA,\mB,\mC,\mD$ are four consecutive intervals.

\EEN
\eeD

\beR  The four sets  $\mA,\mB,\mC,\mD$ correspond to the cases when equality in the HJB equation \eqref{HJB} is attained by at least one of the operators
$-V-P, V'-k,(\mG-q)V,$ and $ 1-V',$ \resp.
Note this generalizes    \cite[Ch. 5.3]{AM14}, where only the last two operators are
considered.
\eeR

\beR One may conjecture that dividend-limited injections
 strategies are of a (recursive) multi-band nature.
 In the case of \expoj, \cite{AGLW} show that the  four sets $\mA,\mB,\mC,\mD$  in the optimal solution  are intervals,  denoted respectively by $(-\I,-a), [-a,0], (0,b),[b,\I)$.

 Cheap equity corresponds the case when $\mC=\emptyset$, and the partition reduces to three sets.
\eeR

\sec{Explicit determination of $a^*,b^*$  when $ F(x):=1- e^{-\mu x}, P > - \fr cq $ \la{s:eqc}}
In this section we turn to the exponential case, where
 explicit formulas for the optimizers $a^*,b^*$ are available. In particular, we will take advantage of properties of the Lambert-W function, which were not exploited in \cite{AGLW}.  Subsequently, in sections \ref{s:ci}, \ref{e:Exp32per} we will show that  exponential approximations
work typically excellently in the general case. Although these results have already been established in [AGLW20], the present formulations have two achievements:
\begin{enumerate}
\item allow an unified formulation of \cite{AGLW} and \cite{Gaj} (via the previously introduced scale functions);
\item make use of a numerical tool (Lambert-W function) to express the optimal quantities of interest $a^*, b^*$.
\end{enumerate}

\ssec{The simplified cost function and  optimality equations}
\beP \la{costexp} {\bf Cost function and  optimality equations in the \expo\ case}
\BEN \im
  \be J_0^{a,b}=
\fr{1-C'(b)\pr{k \; m(a)+P \ovl F(a)} }{(\ovl F(a)) C'(b) + q W_q(b)}=\frac{{\g(b)}-k \;  m(a)-P \ovl F(a)}
{q \th(b)+ \ovl F(a)}, \la{J0}
 \ee
where we put $$ \g(b)=\fr{1}{C'(b)}, \th(b)=\fr{W_q(b)}{C'(b)}.$$

\im Put  \be j(b):= \fr{ \g'(b)}{q  \th'(b)}. \la{jdef}\ee

 For fixed  $a\geq 0$, the optimality equation $\fr{\partial}{\partial b} J_0^{a,b}=0$  \mbw
\be J_0^{a,b}=j(b) \la{j} .\ee

 \im %Put $ J_{0a}:=k a -P$.
 For fixed $k\geq 1$ and $b\geq 0$, at critical points with $a(b)=a^{(k,P)}(b)\neq 0$ \sats\ $\fr{\partial}{\partial a} J_0^{a(b),b}
=0 $ \wmh
 $$ \pp{J_0^{a,b}-(ka -P)}_{a=a(b)}=0.$$ Explicitly,
\be \la{strEqa}
 0=\eta(b,a):=\fr {\g(b)}{\th(b)}- \fr k{\mu \th(b)} F\pr{a}- q\pr{ k a- P}.%=\fr{1}{W_q(b)}-q J_0^{a(b),b} - \fr k{\mu \th(b)}  F\pr{a(b)}.
 \ee

 \im  When $P \geq -\fr c q$ and $ b \geq 0$ is fixed, the solution of \eqr{strEqa} \mbe the principal value of the ``Lambert-W(right)"   function (an inverse of $L(z)=z e^z$) 
 $$[-e^{-1},\I) \ni L_0(z), z \in [-1,\I)$$ \cite{corless1996lambertw,boyd1998global,brito2008euler,pakes2015lambert,vazquez2019psem}
 (this observation is missing in \cite{AGLW}).
\begin{align}\label{Estima*}
&(0,\I) \ni a(b) = \mu^{-1}\pr{-h(b) +L_0\pr{\fr{e^{ h(b)}}{ q \th(b)}}}, \quad
h(b)=h(b,P)=\fr 1{ q \th(b)}-\fr{\mu }{ k}\pr{\fr {  \g(b)}{ q \th(b)}+P }  %>0.
\end{align}

\Itf
\be \la{J0L} J_0^{a(b),b}=\fr k {\mu}\pr{-h(b) +L_0\pr{\fr 1{ q \th(b)}e^{ h(b)}}}-P .\ee

 \im In the special case $b=0$,  \eqref{strEqa} implies that $a=a^{(k,P)}=a^{(k,P)}(0)$ satisfies
the simpler equation
\begin{align}\label{strEqa0}
0=\d_{k,P}(a):=\l \eta(0,a)=\T c -k\pr{a q + \fr{\l}{\mu}(1-e^{-\mu a})  }, \; \T c =c + q P >0,\end{align}
with solution \begin{align} \la{e:L} \mu \; a^{(k,P)}=-g+ L_0\pr{\fr \l q e^{ g}} >0, \; g=h(0)=\fr{\l }{ q}-\fr{\mu }{k q} \T c.
\end{align}

\im At a critical point $(a^*,b^*), a^*>0, b^* > 0$,  \wmh\ both $J_0^{a^*,b^*}  =j(b^*)=k a^* -P\Lra$
\be    a^* =s(b^*), s(b):=\fr{j(b)+P}k, \la{J0a}\ee
and %$b^*$  may be computed from
\be \la{strEqb} %q \th(b) j(b)+ \fr k{\mu} \pr{1- e^{-\mu a}}=\g(b) \Eq
0=\eta(b^*), \quad \eta(b):=\eta(b,s(b))=\fr{\g(b)}{\th(b)} -q  j(b) - \fr k{\mu \th(b)}  F\pr{\fr{j(b)+P}k}=0. \ee

\im The equation $0=\eta(b)$ may be solved explicitly for $P$, yielding
\be P= - \fr k{\mu} \log\pp{1+\fr{q \th(b) j(b) -\g(b)}{\fr k{\mu}}} -j(b).\ee %\fn[4]{This is similar with \eqr{Estima*}, however there $P$ appears also in $h(b,P)$}

\EEN
\eeP

\prf 1. follows from Theorem \ref{p:cost}.1.

2. Let $M(b),N(b)$ denote the numerator an denominator of $J_0^{a,b} := \fr {M(b)}{N(b)} $ in \eqr{J0}. The optimality equation $\fr{\partial}{\partial b} J_0^{a,b}=\fr {N'(b)}{N(b)}\pr{\fr {M'(b)}{N'(b)} -J_0^{a,b}}
=0  $ simplifies to
$$J_0^{a,b} = \fr {M'(b)}{N'(b)}=\frac{{\g'(b)} }
{ q \th'(b)}=j(b).%(=-\fr{C''(b)}{q\pr{ W_q'(b) C'(b)-C''(b)  W_q(b)}}),
$$% and recalling $J_0^{a,b}=k a$ yields the result.

3. \eqr{strEqa} is a consequence  of 1 and of the \sm\ result Corollary 2.

4. See the proof of the particular case 5;  $a\in (0,\I)$ holds since $P \geq - \fr cq \Lra h(b) < \fr 1{ q \th(b)}.$

5. \eqr{strEqa0} follows  from $W_q(0)=\fr 1c,\th(0)=\l^{-1}.$
To get \eqr{e:L},  rewrite the equation \eqr{strEqa0} as $z e^z =\fr \l q e^{ g}, z=\mu a+g$;   $a\in (0,\I)$ holds since $P \geq - \fr cq \Lra g < \fr \l{ q }.$

6. follows from 2. and .3.

7. is \strf.

\qed
\iffalse
\beR $j(b)$ may also be solved explicitly, yielding
\be  {\mu} j(b)=k L_0\pp{
\fr{e^{-P \fr {\mu}k}+ \g(b) \fr {\mu}k -1}{q \th(b)}} .\ee
%From this, $b$ may be solved in principle as well, using the inverse function of $j$.
\eeR
\fi

\beR Note that the \deF\ and \SLG\ solutions $a^*=a(b^*)=\bc 0\\\I \ec$  are always non-optimal, when
$P \geq -\fr c q$ (see \eqr{Estima*}).

\How
 as $k\rightarrow\infty$, $h(b)\rightarrow\frac{1}{q\theta(b)}\neq 0$ and, $a(b)=\mu^{-1}\pr{-h(b)+L_0\pr{h(b)e^{h(b)}}}=0$. This suffices to infer that you get \deF case.
\iffalse
\be {k \to \I} \Lra h(b) \to 0 \Lra  a(b) \to 0 \Lra \bc J_0^{a(b),b} \to -P\\ \eta(b)\to \fr 1{W_q(b)}\pr{J_{0,P}^{DeF}(b)-j(b)}\ec, \for b >0,\ee
\fi

\Oth
\be {P \to \I} \Lra h(b) \to -\I \Lra  a(b) \to \I \Lra \bc J_0^{a(b),b} \to \fr{\g(b)-k/\mu}{q \th(b)}=\fr{1-k C'(b)/\mu}{q W_q(b)}=J_{0,SLG}(b)\\ \eta(b)\to q\pr{J_{0,k}^{SLG}(b)-j(b)}\ec, \for b >0.\ee

Thus, these regimes can be recovered asymptotically. {Let now
$b_k^{*,S},b_P^{*,D}$ denote the unique roots of $\eta(b)=0$ in the two asymptotic cases, which coincide with the classic \SLG\ and \deF barriers.

%, because the respective limits of $\eta(b)$ are proportional to these objectives.

Then, \itm $b^* \leq \min[b_k^{*,S},b_P^{*,D}]$.}\eeR

\iffalse
\beR   The important equation \eqr{J0a} identifies  the optimal buffer associated with a dividends barrier $b$ via the explicit function $s(b)$. In the general framework of \cite{Gaj}, $s(b)$ is only defined implicitly as solution of   \cite[(6)]{Gaj}.

\eeR
\beR Without switching to $\g, \th$, the previous computation is
more complicated
\bea J_0^{a,b}=
\fr{-k \; m(a)C''(b)}{(\ovl F(a)) C''(b) + q W_q'(b)}\eea

\eeR

\fi

\ssec{Existence of the roots of the equations $\eta(b)=0,\d_{k,P}=0$}
The following (new) result discusses the existence of the roots of the equations $\eta(b)=0,\d_{k,P}=0$ introduced in  proposition \eqr{costexp} and relates them  to the  Lambert-W function.

\beP \la{MainP}
\BEN

\im $\th$ increases from
$\th(0)=\fr{1/c}{\l/c}=\fr 1{\l}$ to $\th(\I)=\fr 1{c \Fq - q }$, as we see it in the figure below.%=\fr{\mu \Fq+1}{\l}$ .
  \begin{figure}[H]
    \centering
        \includegraphics[scale=0.7]{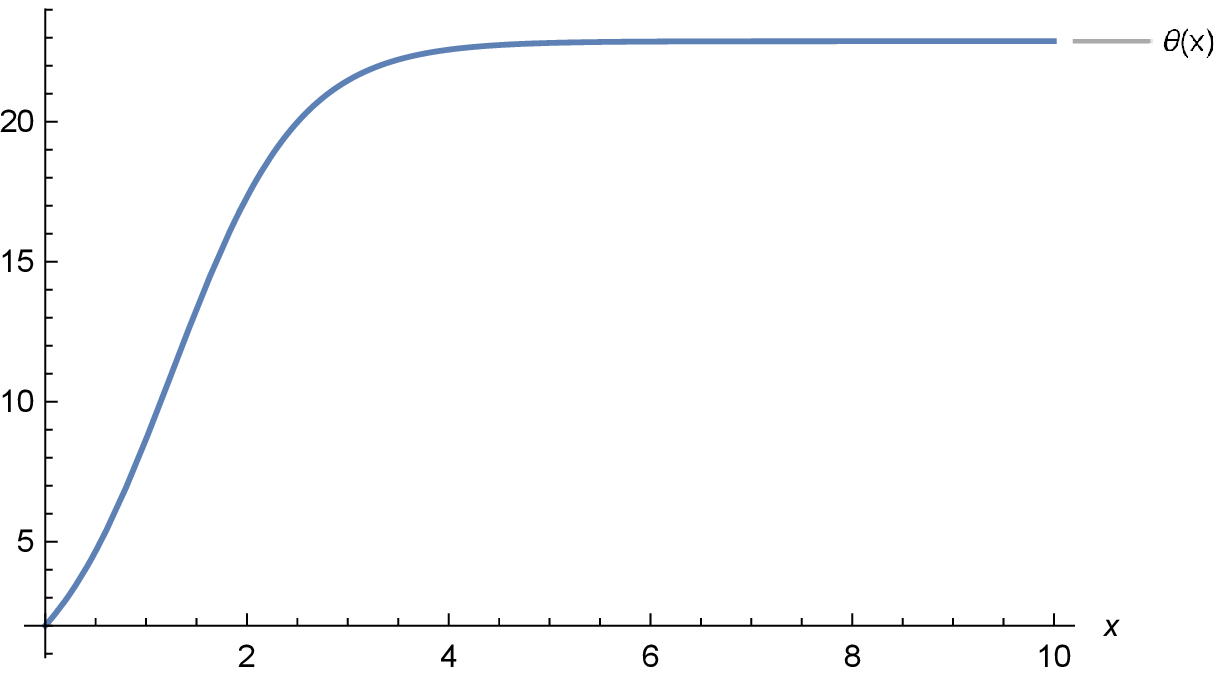}
        \caption{Plot of $\th$ with $\th(0)=2$ and $\th(\I)= 22.8743$, for $\mu=2,\ c=3/4,\ \l=1/2,\ q=1/10$, $ P=1$ and $k= 3/2$.}
        \label{fig:th0}
   \end{figure}
$\g$ is increasing-decreasing (from $\fr c \l$ to $0$), with a maximum at
the unique root of $C''(x)=0$ given by
\be \la{bb}\bar{b}:=\frac{1}{\Phi_q-\rho_-}\log\pr{\frac{\rho_-^2}{\Phi_q^2}},\ee
where $\Phi_q,\rho_-$ denote the positive and negative roots of the \CL\ equation $\k(s)=0$.

The figure below illustrates the plot of the function $\gamma$ and $j({b})$ in which the $\bar{b}$ is represented by the black point.
\begin{figure}[H]
     \centering
        \includegraphics[scale=0.7]{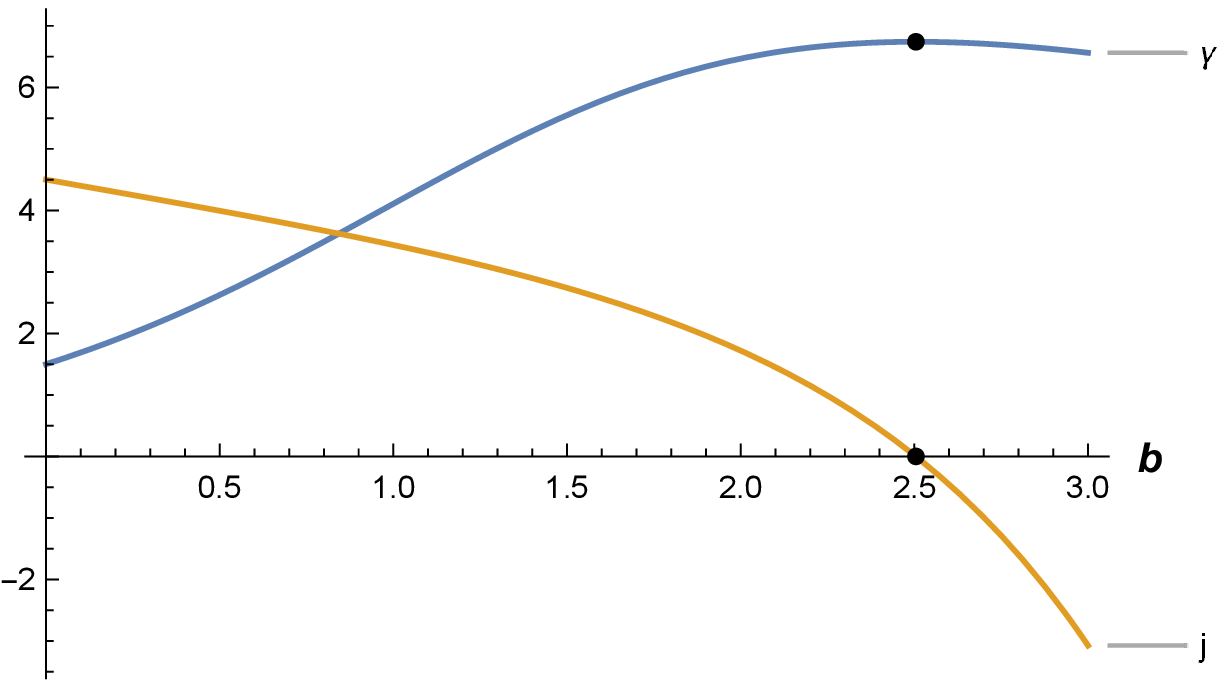}
    \caption{Plots of $j({b})$ and $\gamma({b})$ with $\bar{b}= 2.5046$ and $j(0)=4.5$, for $\mu=2,\ c=3/4,\ \l=1/2,\ q=1/10$, $ P=1$ and $k= 3/2$.}\label{fig:gajb}
\end{figure}

If  $c \mu -  \pr{q +\l} >0$, then $\bar{b}>0$ defined in \eqr{bb} is the unique positive root  of $j({b}) $ and
$\eta\pr{\bar{b}} =\fr{1}{W_q\pr{\bar{b}}}>0$. See the figure \eqr{fig:etbb}

\begin{figure}[H]
    \centering
        \includegraphics[scale=0.7]{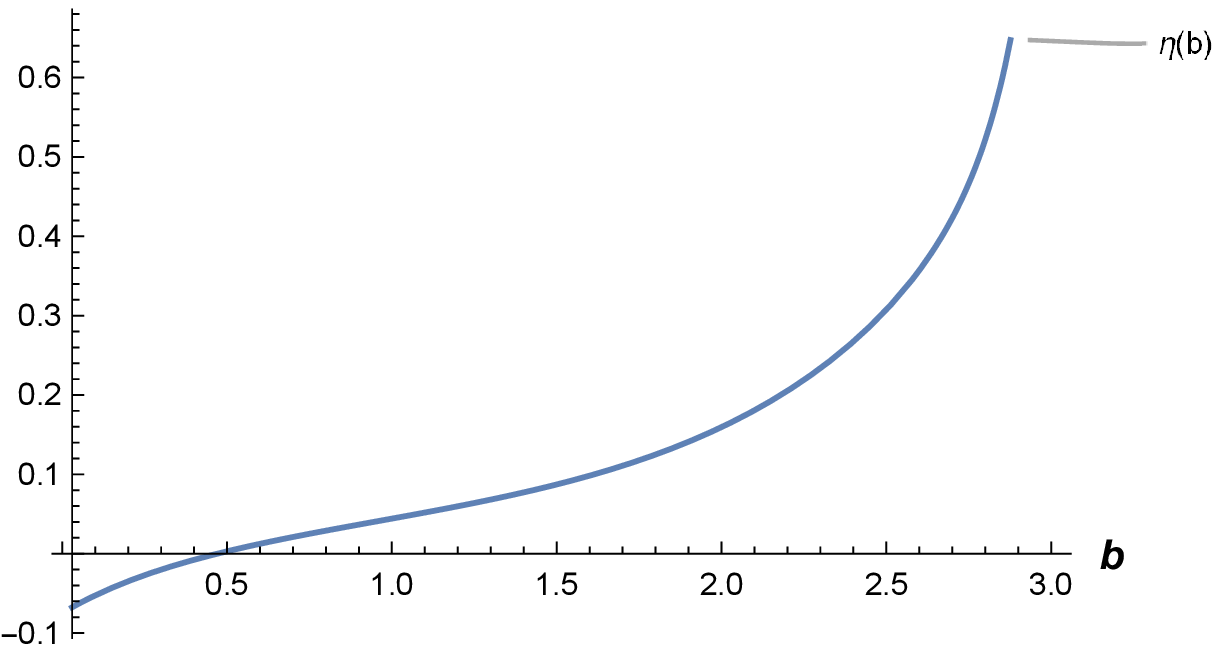}
        \caption{For $\mu=2,\ c=3/4,\ \l=1/2,\ q=1/10$, $ P=1$ and $k= 3/2$, the root of  $\eta({b})=0$ is at $b=0.469843$} %, and $\eta(\bar{b})= 0.337055 >0$ .}
        \label{fig:etbb}
   \end{figure}

The function $j(b)=\frac{\g'(b)}{q \th'(b)}$ is \nne\ and decreasing to $0$ on $[0,\bar b]$, with
\be j(0)=\fr{\l}{\mu q}\fr{-C''(0)}{(C'(0))^2}=
\fr{c \mu -  \pr{q +\l}}{\mu q} \la{j0}. \ee See the figure \eqr{fig:gajb}.

\im Put  \be \la{de} \d_{k,P}:=\d_{k,P}(a^{(k,P)})%=\eta(0) \l
=\d_{k,P}(\fr{j(0)+P}k)=\d_{k,P}(\fr{\T c \mu -  \pr{q +\l}}{k \mu q})=
\fr{\lambda+q-\lambda k\pr{1-e^{-\frac{\T c\mu-\lambda-q}{qk}}}}{\mu},
\ee
and  assume
\be {
\underset{k\rightarrow \I}{\lim}\d_{k,P} = \frac{\l +q}{\mu}- \frac{\l\pr{\mu \T c - (\l+q)}}{q  \mu}= \frac{\pr{\l+q}^2- \l \mu \T c}{q \mu}< 0} \Eq \T c \mu>\lambda^{-1} \pr{\lambda+q}^2. \ee

Then, $\for P> -\fr c q$, the  function $ \d_{k,P}$ is decreasing in $k$ with $\d_{1,P}>0$, and has a unique root
\be \la{kc} k_c=k_c(P):=\fr{q + \l} \l \fr{f}{f + L_0\pr{- f e^{-f}}}>\fr{q + \l} \l,\ee
where
\be \la{f}  f:=\fr{ \l}{q+\l} \fr{\T c\mu-\pr{\lambda+
q }}q> 1 \Eq \T c \mu>\lambda^{-1} \pr{\lambda+q}^2 \Eq   P > P_l:=q^{-1}\pr{\mu^{-1} \lambda^{-1} \pr{\lambda+q}^2 -c}  \ee
(note that the denominator $f+L_0\pr{- f e^{-f}}$ does not equal $0$ since $f >1$ and $L_0$ takes always values bigger than $-1$; or, note that $
-f=L_{-1}\pr{L(-f)}$, where $L_{-1}$ is the other real branch of the Lambert function).

\Fr   %(implying $b^* >0$)
\begin{equation}
\label{Cheapk}\d_{k,P}<0 \Eq k >k_c(P).\end{equation}

\im \Itf
  $\eta(b)=0$  has at least one solution of in $(0,\bar{b}]$ iff
\be \la{de} \eta(0)=\fr c \l- \fr 1 \l \pr{c - \fr {q +\l}\mu}-\fr {k}\mu F\pr{a^{(k,P)}}=\fr {1}{\l \mu} \pr{\l+  {q}-\l k F\pr{a^{(k,P)}}}=\fr {1}{\l }\d_{k,P}<0 \Eq k >k_c.
\ee

The first such solution will be denoted by $b^*$. % Explicitly,  the structure equation \eqr{strEqb} may be rewritten as:\bea (\fr {\mu } k j(b)+ h) e^{\fr {\mu } k j(b)+ h}=\fr {e^{h }}{q \th(b)}  \Lra \fr {\mu } k (j(b)+P) = -h +L_0(\fr {e^{h }}{q \th(b)} ).\eea
%\be     \eta(b):=\fr {1}{W_q(b)}\pr{1 - k\fr {C'(b)}{\mu } F\pr{s(b)}}   -k  q s(b)+ =0. \la{strEqb}\ee

\EEN

\eeP
\prf
 For 1.  see  \cite[Proof of Theorem 11, A2]{AGLW}.

  2. By using the assumption $ \T c \mu > \l^{-1} (\l+ q)^2$ we get $\T c \mu \geq \l +q$, and $k \in [1,\I) \rightarrow \d_{k,P}$ is decreasing.

  Put $d =\fr{\T c\mu-\pr{\lambda+
q }}q$. The inequality $\d_{k,P}<0$ (see \eqr{de}) \mbr

\bea   e^{-\frac{d}{k}} <  1 -\fr{q+\l}{\l k} \Eq 1 < e^{\frac{d}{k}}
\pr {1 -\fr{(q+\l)}{d \l } \fr d k}:=e^{z} \pr {1 -z/f}.
\eea

Rewriting the latter as $-f>e^{z} \pr {z -f}$  we recognize, by putting $z=y+f$,  an inequality
reducible   to $ y e^{y} <- f e^{-f}. $ The   solution is $$y <L_0\pr{- f e^{-f}},$$ where $L_0$ is the principal branch of the Lambert-W function.
%\red{ Note that the argument is negative, and so a branch must be chosen}.

The final solution
is \eqr{Cheapk}, where we may
note that the  variables $k,P$ have been separated.

 3. is \strf.

\qed

 \beR

The function $[1,\I) \ni f \mapsto  \fr{f}{f + L_0\pr{- f e^{-f}}} \in (1,\I)$ blows up at $f=1$, and converges to $1$ when $f \to \I$
  (or when either $\mu $ or $\T c=c+  q P$ are large enough) as may be noticed in the figure below, which blows up at the value $P_l:=-4/5 .$ Note also that when  $f$ (or one of $c,P,\mu$ are large enough),  $k_c$ given by \eqr{kc} stabilizes to the equilibrium $\frac
{\lambda+q}{\lambda}=6/5$;  this is related to  \cite[Lemma 2]{APP}, \cite[Lemma 7]{kulenko2008optimal}, who obtain the same condition for $b^*=0$ (without buffering capital injections). Intuitively, under these conditions,  buffering is not crucial.

At the other end, as $f $ tends to its lower limit and to the regime B,   the notion of equity expensiveness  vanishes, and  $k_c \to \infty$.
\figu{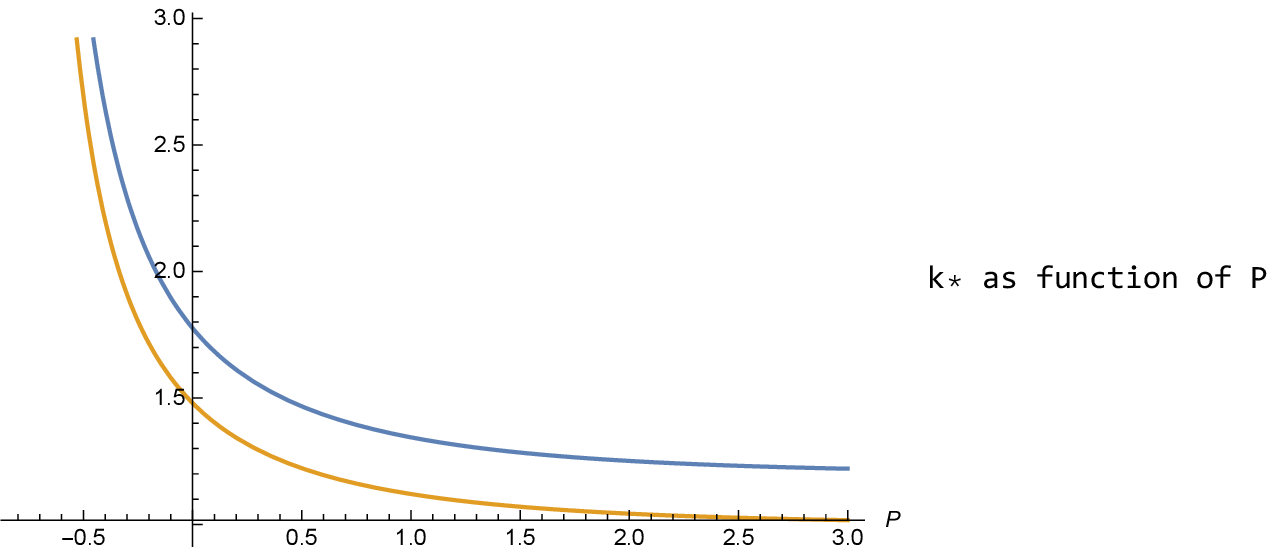}
  {$k_c$ as function of $P$, for several values of $c$, with the vertical asymptote at $P_l$ fixed.}{}

\eeR

%\newpage
The next  two figures illustrate how $k_c$ blows up at the critical values $q_l :=(1/2) (-2 \l  + P \l  \mu  +
   \sqrt\l  \sqrt\mu  \sqrt
    4 c - 4 P \l  + P^2 \l  \mu )$  and $\l_l :=  \left(c \mu -\sqrt{(-c \mu +\mu  (-P) q+2 q)^2-4 q^2}+\mu  P q-2 q\right)$
(represented by  red points in the figures below). The dark (blue) parts correspond to the regime $A$.

\begin{figure}[H]
    \centering
    \begin{subfigure}[b]{0.45\textwidth}
        \includegraphics[width=\textwidth]{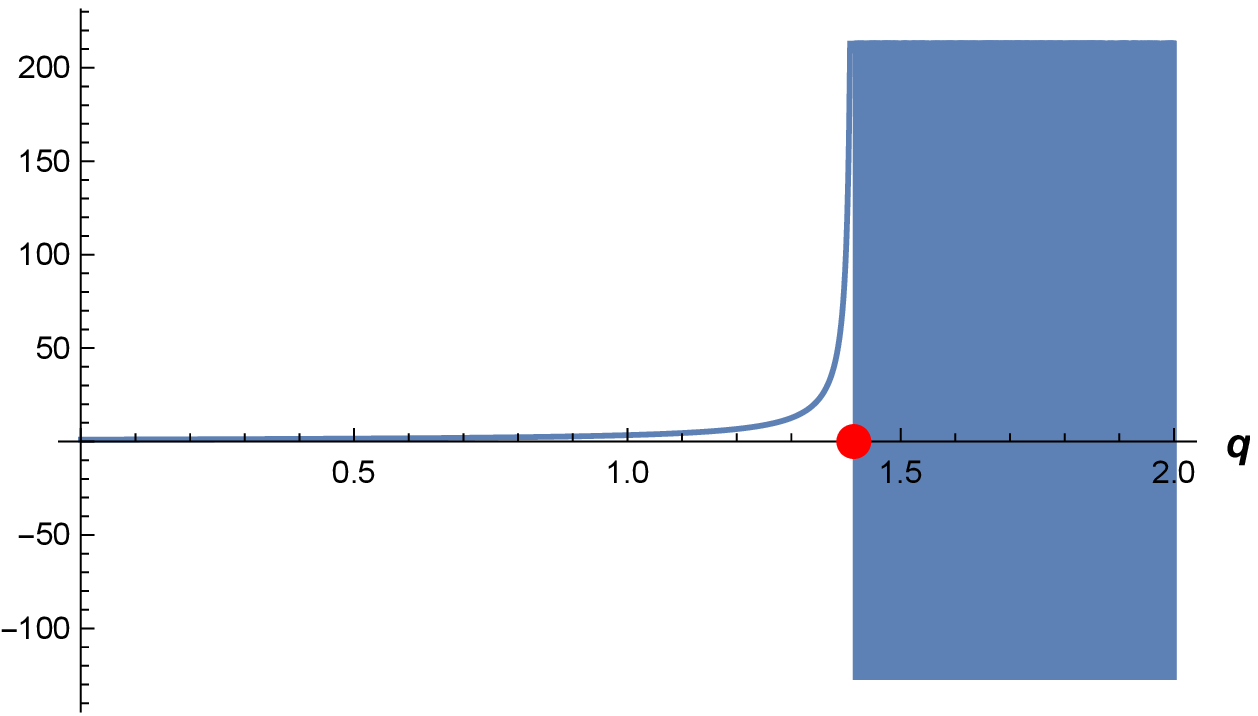}
        \caption{$k_c$ defined in \eqr{kc}as function of $q$, for $\mu=2,\ c=3/2,\ \l=1$ and $ P=1$; $q_l= \sqrt{2}$.}
        \label{fig:kcq}
    \end{subfigure}
    ~ %add desired spacing between images, e. g. ~, \quad, \qquad, \hfill etc. 
    \begin{subfigure}[b]{0.45\textwidth}
        \includegraphics[width=\textwidth]{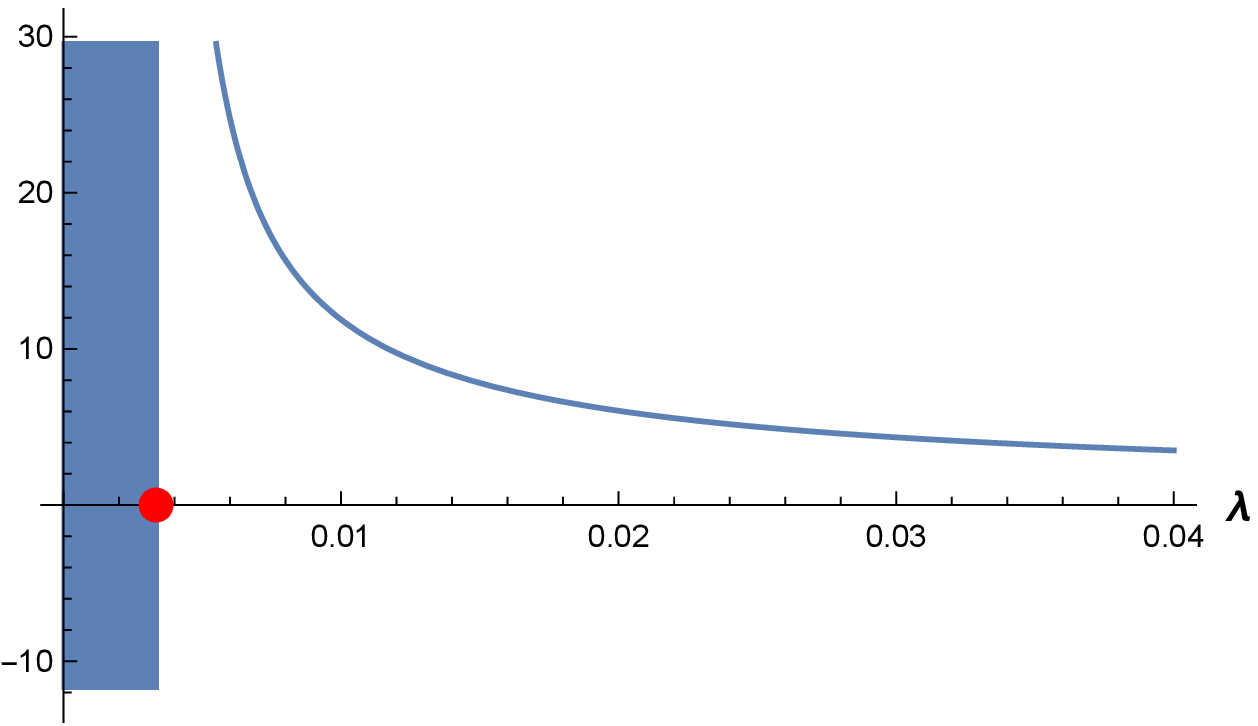}
        \caption{$k_c$ as function of $\l$, for $\mu=2,\ c=3/2,\ q=1/10$ and $P=1$.}
        \label{fig:ptl}
    \end{subfigure}
    \caption{$k_c$ as a function of $q$, $\lambda$}\label{fig:kc}
\end{figure}

\section{Which exponential   approximation? \la{s:DV}}
\ssec{Three de Vylder-type exponential approximations for the \rp\ \la{s:DVr} }
In the simplest case of \expoj\  of rate $\mu$ and $\s=0$, the  formula for the \rp\ is
\begin{equation} \la{Ruiex}
\Rui(x) = P_x[ \exists t \geq 0: X_t < 0]= \frac{1}{1+  {\th}} \ exp \left( - \frac{x  {\th} \mu }{1+  {\th}} \right)= \frac{1}{1+  {\th}} \ exp \left( - \frac{x  {\th} m_1^{-1}}{1+  {\th}} \right),
\end{equation}
where $\th=\fr{c -\l m_1}{\l m_1}$ is the \lc.
By plugging the correct mean of the claims in the second formula yields the {\bf simplest approximation} for \procs\ with finite mean claims.

More sophisticated is the {\bf Renyi exponential
approximation}
\begin{equation} \la{Re}
\Rui_R(x) =  \frac{1}{1+  {\th}} \ exp \left( - \frac{x  {\th} \H m_2  ^{-1}}{1+  {\th}} \right), \H m_2  =\frac{m_{2}}{2m_{1}};
\end{equation}

This formula can  be obtained as a \tPd\   of the \LT, which conserves also the value $\Rui(0)=(1+\th)^{-1}$ \cite{AP14}. It may be also derived heuristically from the first formula in \eqr {Ruiex}, via replacing $\mu$ by the correct ``excess mean" of the  {\bf excess/severity density} $$f_e(x)=\fr{\ovl F(x)}{m_1}=\fr{1- F(x)}{m_1},$$
which is known to be $\H m_2  $. Heuristically, it makes more sense to approximate $f_e(x)$  instead of the original density $f(x)$, since $f_e(x)$ is a monotone function, and also an important component of the \PK formula for the \LT\ $\H \Rui(s)=\int_0^\I e^{- s x} F( dx)$ -- see \cite{ramsay1992practical,AP14}.

More moments are put to work in the {\bf \deV}
\begin{equation}
\Rui_{DV}(x) = \frac{1}{1+ \T {\th}} \ exp \left( - \frac{x \T {\th} \H m_3 ^{-1}}{1+ \T {\th}} \right), \; \H m_{3}:=\frac{m_{3}}{3m_{2}}, \;  \T {\lambda}
= \frac{9 m_{2}^{3}}{2 m_{3}^{2}}\lambda, \; \;  \T c= c -\l m_1+ \T \l \H m_3 , \;  \T \th =\frac{2 m_{1}m_{3}}{3m_{2}^{2}}\th=\frac{\H m_{3}}{\H m_{2}}\th.%= c -\l m_1+  \l \frac{ m_{2}^{2}}{2 m_{3}}.
\end{equation}
Interestingly, the result may be expressed in terms of the so-called "normalized moments" \be {\H m_i}=\fr{m_i}{i \; m_{i-1}}\ee  introduced in \cite{bobbio2005matching}.

The \deV\ parameters above may be obtained either from  \BEN \im  equating the first three cumulants of our process to those of a process  with  exponentially distributed claim sizes of mean $\H m_3 $,  and  modified $\l, c$  \cite{de1978practical}  (however   $p=c -\l m_1=E_0[X_1]$ must be conserved, since this is the first cumulant), or
 \im a Pad\'e approximation  of the Laplace transform  of the \rps\ \cite{avram2011moments}.
\EEN

The  second derivation via Pad\'e shows that higher order approximations may be easily obtained as well. They might not be admissible, due to negative values, but packages for ``repairing" the non-admissibility are available -- see \fe  \cite{dumitrescu2016modeling}.

The first derivation of the  de Vylder approximation is a {\bf process approximation} (i.e., independent of the problem considered); as such, it may  be  applied to other functionals of interest  besides \rps ($W_q(x)$, dividend barriers, etc), simply by plugging the modified \para s in the exact formula for the \rp\ of the simpler process.

\iffalse
\section{Can we minimize the \rp, by proportional \rei? \la{s:DVr} }

Let us consider first the mean premium principle, with $\th_R \geq \th$ and retention level (potentially dependent on $x$). As \wk, we may forget about the reinsurer, and say instead that we modify the premium to
and the \lc\ to
\be \la{Tc} \T c=\l m_1 \pr{1+ \th - (1-a) \pr{ \th_R +1}}=\l m_1 \pr{1+ \T \th},\ee
\be \la{pos} \T \th=\th - (1-a) \pr{ \th_R +1},\ee
with $a \geq \fr {1 +\th_R-\th}{1+ \th_R } $, to ensure positivity in \eqr{pos}.
 \Fr we must also multiply the mean by $a$. %the level of the natural exponential  approximation

 The \rp\ has one critical point with $a$ satisfying  a quadratic equation (but, taking $a$ at the critical point  yields negative $\T \th$). A careful analysis reveals that the optimum is at $a=1$?.

 This continues to be the case for the mean-variance premium principle (Shihao?).
 The optimization becomes non-trivial for excess of loss \rei?.
\fi

\ssec{Three two point Pad\'e approximations of the \LT\  $\H W_q$ of scale function \la{s:DVW}}

The simplest approximations for the scale function  $W_q(x)$ will now be derived heuristically from the following example.

 \beXa
  {\bf The Cram\'{e}r-Lundberg model with exponential jumps \la{s:exp}}
Consider  the Cram\'{e}r-Lundberg model
 with exponential jump sizes with mean $1/\mu$, jump
rate $\lambda$, premium rate $c>0$,
and Laplace exponent
$\k(s)=s \le(c-\fr{\lambda}{\mu+s}\ri)$. Solving $\k(s)-\q=0 \Eq c s^2 + s(c \mu -\l -\q) - \q \mu=0$ for $s$ yields two distinct solutions $\g_{2} \leq 0 \leq \g_{1}=\Fq$ given by
\begin{align*}
\g_{1}(\mu,\l,c)=\g_1 =& \fr{1}{2c} \left(- \left(\mu c -\lambda - \q\right) + \sqrt{\left(\mu c -\lambda - \q \right)^2 + 4\mu \q c} \right),\\
\g_{2}(\mu,\l,c)=\g_2 =& \fr{1}{2c} \left(- \left(\mu c -\lambda - \q\right) - \sqrt{\left(\mu c -\lambda - \q \right)^2 + 4\mu \q c} \right).
\end{align*}

The $W$ scale function is:
\be \la{Wexp}  W_{\q}(x) = \fr{A_1 e^{\g_{1}x} - A_2 e^{\g_{2}x}}{c(\g_{1}-\g_{2})}  \Eq \H W_{\q}(s) = \fr{s+ \mu}{c s^2 + s(c \mu -\l -\q) - \q \mu},\ee
where $A_{1} = \mu + \g_{1}, A_{2} = \mu + \g_{2}$.

Furthermore, it is \wkt    the function $ W_{\q}'(x)$ is in this case
unimodal  with global minimum at
\be \la{ob} b_{DeF} = \frac{1}{\g_{1} - \g_{2}}
\begin{cases}\log
\frac{(\g_{2})^2 A_{2}}{(\g_{1})^2 A_{1}}=\log
\frac{(\g_{2})^2(\mu +\g_{2})}{(\g_{1})^2(\mu +\g_{1})} \quad &\text{if $ W_{\q}''(0) <  0 \Eq (\q+\lambda)^2$}-
c\lambda\mu < 0\\ 0 & \text{if $ W_{\q}''(0)\geq 0 \Eq
(\q+\lambda)^2- c\lambda\mu \geq 0$}\end{cases}, \ee since
$ W_{\q}''(0) =
\fr{(\g_{1})^2(\mu +\g_{1})-(\g_{2})^2(\mu +\g_{2})}{c(\g_{1}-\g_{2})}
=\fr{ (\q+\lambda)^2- c\lambda\mu}{c^3}$ and that the optimal strategy for the \deF \prob \ is
 the barrier strategy at level $b_{DeF}$ (see \fe \cite{APP}, \cite[Sec. 3]{AGV}).
\iffalse
 Also, the optimal barriers in the presence of a final penalty $P$, and of reflection with proportional costs $k$, \saty \resp
 \be \bc P  \q \Delta_\q^{(W)}(b_P)=-W_q''(b_P)\Eq P  \q \mu \Tl c^{-2}
     e^{(\Tl + \Tq -\mu) b_P}=-W_q''(b_P) \\k \Delta_\q^{(ZW)}(b_k)=W_q'(b_k) \Eq k \Tl c^{-2}
     e^{(\Tl + \Tq -\mu) b_k}=W_q'(b_k)
  \ec. \ee
\fi
\eeXa

 Plugging now the respective  \paras\ of the de Vylder type approximations in  the exact formula \eqr{Wexp}  for  the \CLp\ with \expoc,  we obtain  three approximations for $\H W_q$:
 \BEN \im   ``Naive exponential" approximation obtained by plugging  $\mu^{-1} \rightarrow m_1$ in \eqr{Wexp}
    (as was done, for a different purpose) in \eqr{Ruiex}
    \im  Renyi\fn[4]{This is called  DeVylder B)  method in \cite[(5.6-5.7)]{GSS}, since it is the  result   of  fitting the first two cumulants of the risk process.}, obtained by plugging $\mu^{-1} \rightarrow \H m_2, \l_R \rightarrow \l \fr{ m_1}{\H m_2}$ (since $c$ is unchanged, the latter equation is equivalent to the conservation of $\rho=\fr{\l m_1}c,$ and to the conservation of $\th$, so this coincides with the  Renyi ruin approximation used in \eqr{Re}.)
  \im   De Vylder, obtained by plugging $\mu^{-1} \rightarrow \H m_3, \T \l \rightarrow \l \fr{ 9 m_2^3}{2 m_3^2}, \T c=c- \l m_1+ \T \l \H m_3$. \EEN

  \beR  In the case of \expoc,  these three  approximations are exact, by definition
(or check that for \expoc \ all the normalized moments are equal to $\mu^{-1}$). \eeR

 \beR The conditions for the non-negativity of the barrier is $W_{\q}'' (0_+) <0 \Eq (\frac {\lambda+ \q }  c )^2 < \frac {\lambda}  c f(0)$. Here, this condition is \satd\ for   the exact when $\th > \frac{(\lambda +q)^2 (1- \rho)}{\lambda^2 f(0) m_1}$. %For the approximations, it is \satd\ at ...

\eeR

It is shown in \cite[Prop. 1]{AHPS} that the three de Vylder type approximations   are  two-point Pad\'e approximations of the \LT\  (hence   higher order generalizations are immediately available).

 We recall that  two-point Pad\'e approximations incorporate into the \Pd\  two initial values of the function (which can be
derived easily via the initial value theorem, from the \PK \LT):
\beq \la{W0} &&W_q(0_+)=  \lim_{s \to \I}s \H W_q(s)= \frac 1 {\c}, \\&&
W_q'(0_+)= \lim_{s \to \I}s\left( \fr{s}{\k(s) -\q}- W_{\q} (0_+) \right)=   \frac {\q + \l} {\c^2}. \la{W0p}\eeq

In our case, incorporating both $W_q(0_+),
W_q'(0_+)$ leads to the natural exponential  approximation which is therefore the best near $x=0$.
Incorporating none of them yields the \deV, which is the best asymptotically. Incorporating only $W_q(0_+)$ leads to Renyi, which is expected to be the best
in an intermediate regime.

\Nt when  the jump distribution has a density  $f$, \ith:\fn[6]{This equation is important in establishing the nonnegativity of the optimal dividends barrier.}
\begin{equation} \label{e:secder}
\begin{aligned}
 W_{\q}'' (0_+) &=\lim_{s \to \I} s \left( s \left(\fr{s}{\k(s) -\q}- W_{\q} (0_+)\right) -  W_{\q}' (0_+) \right)=    \fr 1 c \Big( (\fr {\lambda+ \q }  c )^2 - \fr {\lambda}  c f(0) \Big).
\end{aligned}
\end{equation}

Thus,  $W_q''(0)$ already requires knowing $f_C(0)$ (which is a rather delicate task starting from real data); \thr\ we will not
incorporate into the \Pd\  more than two initial values of the function.

%\iffalse

 We  recall below in Proposition \ref{p:AH}  three types of two-point \Pds\ \cite[Prop. 1]{AHPS}, and particularize them to the case when the denominator degree is $n=2$ (which are  further illustrated below).

\beP \la{p:AH}  {\bf Three  {matrix exponential} approximations for the scale function}.
\BEN \im To secure both the values of $W_q(0)$ and $W_q'(0)$, take into account \eqr{W0} and \eqr{W0p}, i.e. use the \Pd  $$\H W_q(s) \sim \fr{\sum_{i=0}^{n-1} a_i s^i}{c s^n+ \sum_{i=0}^{n-1} b_i s^i}, a_{n-1}=1, b_{n-1}=c a_{n-2}-\l -\q.$$
For $n=2$ we recover the ``natural exponential " approximation of plugging  $\mu \to \fr 1{m_1}$ in \eqr{Wexp}:
\be \la{W2zz} \H W_q(s) \sim
   \frac{\frac{1}{m_1}+s}{c s^2+ s \left(\frac{c}{m_1}-\lambda
   -q\right)-\frac{q}{m_1}},\ee
   used also (for a different purpose) in \eqr{Ruiex}.

   \im To ensure only $W_q(0)=\fr 1 c$, we must   use the \Pd  $$\H W_q(s) \sim \fr{\sum_{i=0}^{n-1} a_i s^i}{c s^n+ \sum_{i=0}^{n-1} b_i s^i}, a_{n-1}=1.$$
For $n=2$, we find
\be \la{W2z} \H W_q(s) \sim\frac{\frac{2 m_1}{m_2}+s}{c s^2 +\frac{s \left(2 c m_1-2
   \lambda  m_1^2-m_2 q\right)}{m_2}-\frac{2 m_1
   q}{m_2}}=\frac{\frac 1 {\H m_2  }+s}{c s^2 +{s \left(\fr c {\H m_2  }-
   \lambda \fr{ m_1}{\H m_2  }- q\right)}-\frac{
   q}{\H m_2  }},\ee
   where  $\H m_2  =\H m_2= \frac {m_2}{2 m_1}$ is the first moment of the excess density $f_e(x)$.
      Note that it equals the \sf\ of a \proc\
    with  \expoc\ of rate $\H m_2  ^{-1}$ and with    $\l$ modified to  { $\l_R=\lambda \fr{ m_1}{\H m_2  }$.  Since $c$ is unchanged, the latter equation is equivalent to the conservation of $\rho=\fr{\l m_1}c,$ and to the conservation of $\th$, so this coincides with the  Renyi approximation\fn[4]{This is called  DeVylder B)  method in \cite[(5.6-5.7)]{GSS}, since it is the  result   of  fitting the first two cumulants of the risk process.} used in \eqr{Re}.}

 \im The pure \Pd\  yields for $n=2$
\beq \la{W2} && \H W_q(s) \sim \frac{ s+\fr{3 m_2}{m_3}}{s^2 \left( c- \l m_1 + \lambda  \fr {3 m_2^2}{2
   m_3}\right)+s \left( c\fr{3 m_2}{m_3}-\fr{3 m_1 m_2}{m_3} \lambda  - q\right)-\fr{3 m_2}{m_3} q}\no \\&&=\frac{ s+\fr{1}{\H m_3}}{ \T c s^2
+s \left( \T c \fr{1}{\H m_3}-\T \l  - q\right)-\fr{1}{\H m_3} q}, \; \T c= c -\l m_1 + \T \l    \H m_3, \T \l= \l \fr{9 m_2^3}{2 m_3^2}. \no
%\\&&=  \frac{ s+\fr{1}{\H m_3}}{s^2 \left( c -\l m_1 + \T \l_L    \H m_2 \right)
  % +s \left( c \fr{1}{\H m_3}-\Tl_L  - q\right)-\fr{1}{\H m_3} q}, \; \Tl_L=\fr{3 m_1 m_2}{m_3} \lambda=\fr{ m_1 }{\H m_3} \lambda.
  \eeq

   {Note that both the coefficient of $s^2$ in the denominator coincides with the coefficient $\T c$ in the  classic \deV, since $\T \l \H m_3=  \l \fr{9 m_2^3}{2 m_3^2} \fr{m_3}{3 m_2}= \l \fr {3 m_2^2}{2 m_3}$,   and so does  the coefficient of $s$, since
$$ c\fr{3 m_2}{m_3}-\fr{3 m_1 m_2}{m_3} \lambda  = \T c \fr{1}{\H m_3}-\T \l =
 \pr{c -\l m_1 + \T \l    \H m_3} \fr{1}{\H m_3} - \T \l. $$}

\EEN
\eeP

%\fi

\section{Examples of computations involving scale function and dividend value approximations\la{s:ci}}

Our goal in this section is to investigate whether exponential approximations
 are precise enough to yield reasonable estimates for quantities important in control like \BEN \im  the dominant exponent $\Phi_q$ of $W_q(x)$   \im the last local minimum of $W_q'(x)$,  $b_{DeF}$, which yields, when being the global minimum, the optimal De Finetti barrier \im  $ W''_\q(0) $, which determines if  $b_{DeF}=0$ \im the functional $J_0$ yielding the maximum \div s with \ci.\EEN  %Note that other \perf\ measures like the dual optimal dividends barrier, and the reflected optimal dividends barrier could be investigated as well.

 All the examples considered involve a \CL\ model with rational \LT\ $\H{W}_q(s)$ (since in this case, the computation of $W_q,Z_q$ is fast and in principle arbitrarily large precision may be achieved with symbolic algebra systems).

\BEN \im   For the first three \prob s,  we will use de Vylder type approximations.
Graphs of $W'_q$, $W''_q$  and  some tables summarizing the simulation results will be presented. We note that in most of the cases that we observed, the \deV\ of $\Fq$ deviates from the exact value the least -- see \fe Table 2. For the De Finetti barrier, the  "winner" depends on the size of $b_{DeF}$. Unsurprisingly, when near $0$, the natural exponential approximation  wins, and as $b_{DeF}$ increases, Renyi and subsequently the \deV\ take the upper hand -- see \fe Table 3.

\im For the computation of $J_0$, we provide, besides the exact value, also two approximations:
\begin{enumerate}
\item For a given density of claims $f$ one computes an exponential density approximation $f_e(x)=\frac{1}{m_1}\text{exp}(-\frac{x}{m_1})$ where $m_1$ is the first moment of $f$. Subsequently,   $W$, $Z$,  $J_0$ and $a,b$ are obtained using  the exponential \app\ $f_e$. Quantities obtained by this method would be referred to with an affix `expo pure'.
\item For a given density of claims $f$, the value function is computed via the formula which assumes exponential claims in equation \ref{costexp}, but the
    "ingredients" $W$, $Z, \ovl F$ and the mean function $m$  are the correct ingredients  corresponding to our original density $f$. Quantities obtained by this method would be referred to with an affix `expo CI'.
\end{enumerate}

It turns out that the pure expo approximation works better for large $\th$, and the
correct ingredients approximation works better for small $\th$.

Note that we only included tables illustrating  approximating $J_0$  for the first two  examples, to keep the length of the paper under control, but similar results were obtained for the other examples.
\EEN
\iffalse

 \im Use  the ``exact exponential objective function" $J_0^{s(b),b}$ defined   replacing $a$ in \eqr{J0} by $s(b)$, but plug in the non exponential
ingredients $W_q,C_q,\ovl F, m(a)$ of our non-exponential examples.   Then,  optimize this numerically. This is our best approximation, but it is totally "ad-hoc"\fn[4]{We may interpret this approximation  as an attempt to complete a scale based approximation before reaching ruin
 by  taking into account also a better estimate of the severity of ruin/overshoot distribution.}.

 \im  the natural exponential , Renyi and  classic \deV: write a script ``WZforEx" (alias ``Wex") which computes the scale functions $W_q,C_q$ and $\eta(b)$ \cite{str}, \cite{et0}, and apply it with $\mu$ replaced by appropriate parameters corresponding to the natural exponential , Renyi and  classic \deV. This approximations work less
 well, but have clearer interpretations, at least the  natural exponential  (approximate the claims) and the  Renyi \app\ (approximate the equilibrium distribution of the claims).

 We will first compare the optimal values with that of the (exact) numeric optimization  of the matrix-exponential objective function \eqr{J0PH}, and provide relative errors.

We will also compare our three approximations to the previous exact results for the optimal \deF and \SLG\ value functions.\fn[4]{For a reference on the \deF and \SLG\ objectives (which only optimize  the dividend barrier $b$), see \fe \cite[Sec. 2]{AGR}.}

\Ito the best approximation is the classic \deV.

\fi

\subsection{A Cram\'{e}r-Lundberg process with hyperexponential claims of order 2} \label{e:MixExp83}
We take a look at a Cram\'{e}r-Lundberg process with density function
$f(x)=\frac{2}{3}e^{-x}+\frac{2}{3} e^{-2 x}$ with $\l=1$, $\th=1$ and $q=\frac{1}{10}$.

\begin{figure}[H]
    \centering
    \begin{subfigure}[b]{0.45\textwidth}
        \includegraphics[width=\textwidth]{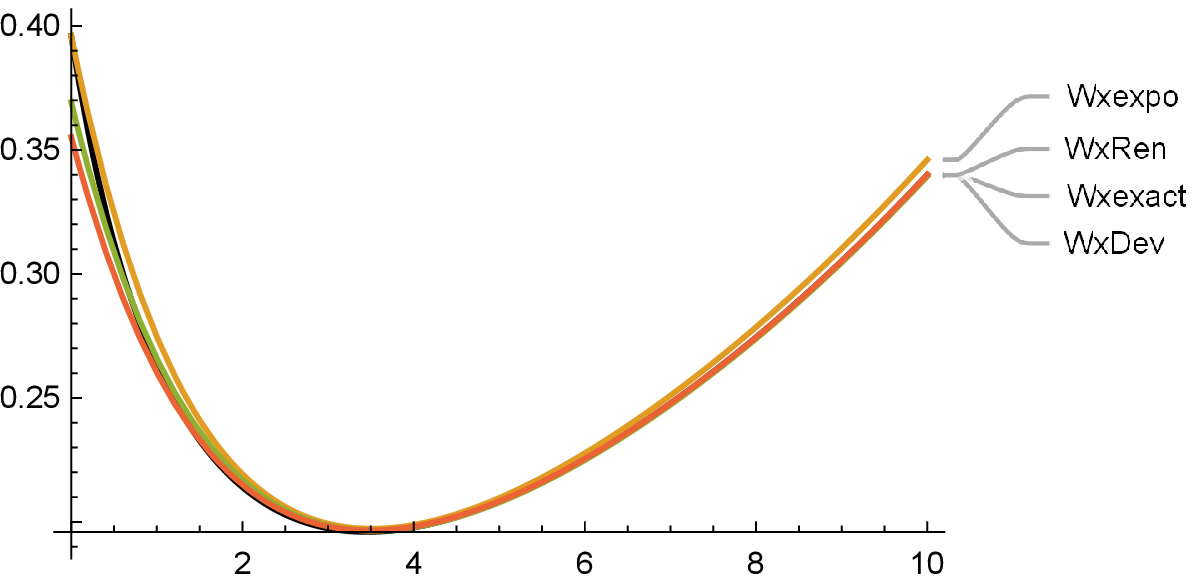}
        \caption{$W'_q(x)$}
        \label{fig:MixExp12W1}
    \end{subfigure}
    ~
    \begin{subfigure}[b]{0.45\textwidth}
        \includegraphics[width=\textwidth]{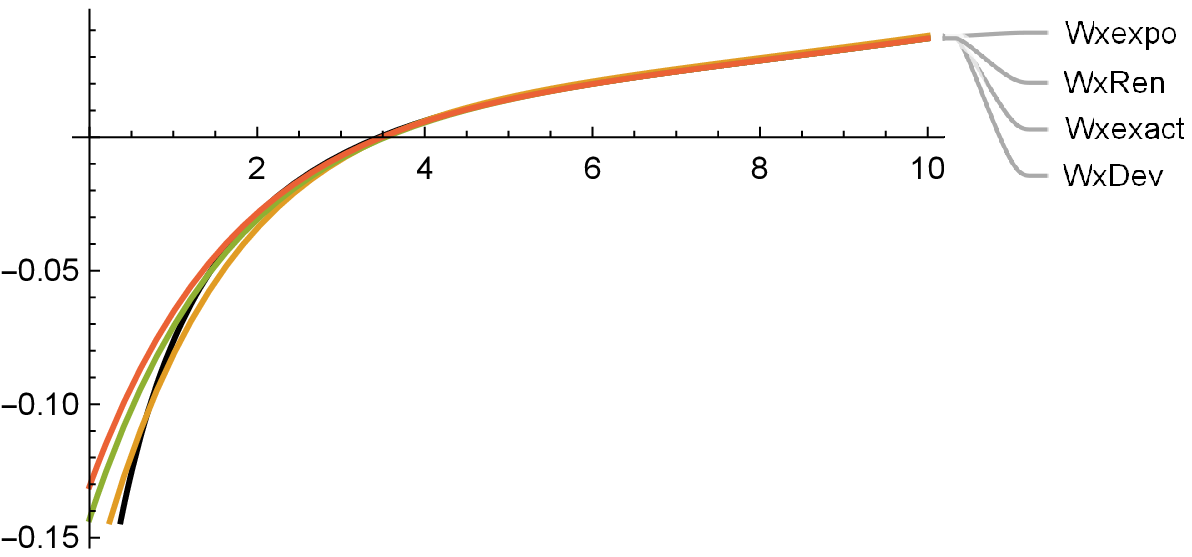}
        \caption{$W''_q(x)$}
        \label{fig:MixExp12W2}
    \end{subfigure}
    \caption{Exact  and  approximate plots of $W'_q(x)$ and $W''_q(x)$  for $f(x)=\frac{2}{3}e^{-x}+\frac{2}{3} e^{-2 x}$, $\theta=1$, $q=\fr{1}{10}$.}\label{fig:MixExp12}
\end{figure}

\begin{table}[H]
\begin{tabular}{|l|l|l|l|l|}
\hline
       & \begin{tabular}[c]{@{}l@{}}Dominant   exponent \\ $\Phi_q$\end{tabular} & \begin{tabular}[c]{@{}l@{}}Percent   relative error\\ ($\Phi_q$)\end{tabular} & \begin{tabular}[c]{@{}l@{}}Optimal barrier\\ $b_{DeF}$\end{tabular} & \begin{tabular}[c]{@{}l@{}}Percent   relative error\\ ($b_{DeF}$)\end{tabular} \\ \hline
Exact & 0.110113                   & 0                                 & 3.45398              & 0                             \\ \hline
Expo  & 0.110657                   & 0.494313                          & 3.51173              & 1.67191                       \\ \hline
Dev   & 0.110115                   & 0.00195933                        & 3.48756              & 0.972251                      \\ \hline
Renyi & 0.110078                   & 0.0321413                         & 3.5323               & 2.26744                       \\ \hline
\end{tabular}
\caption{Exact and approximate values of $\Phi_q$ and $b_{DeF}$ for $f(x)=\frac{2}{3}e^{-x}+\frac{2}{3} e^{-2 x}$, $theta=1$, $q=\fr{1}{10}$, as well as percent relative errors, computed as the absolute value of the difference between the approximation and the exact, divided by the exact, times 100.
Relative errors for $\Phi_q$ are less than $0.5\%$, with the pure exponential approximation proving to be the worst and the DeVylder the best approximations, respectively. The optimal barrier $b_{DeF}$ is also best approximated by DeVylder, with Renyi being the worst at $2.26\%$ .}
\label{table:MixExp12}
\end{table}

\begin{table}[H]
\begin{tabular}{|l|l|l|l|l|}
\hline
$\theta$ & Closest approximation & $\Phi_q$   exact & $\Phi_q$ approximation & \% error   $\Phi_q$ \\ \hline
1          & Dev                     & 0.110113     & 0.110115 & 0.00195933      \\ \hline
0.9        & Dev                     & 0.120328     & 0.120331 & 0.00269878      \\ \hline
0.8        & Dev                     & 0.132452     & 0.132457 & 0.00380056      \\ \hline
0.7        & Dev                     & 0.147017     & 0.147025 & 0.00548411      \\ \hline
0.6        & Dev                     & 0.16475      & 0.164763 & 0.00812643      \\ \hline
0.5        & Dev                     & 0.186652     & 0.186675 & 0.0123901       \\ \hline
0.4        & Dev                     & 0.214122     & 0.214163 & 0.0194631       \\ \hline
0.3        & Dev                     & 0.249118     & 0.249196 & 0.0315039       \\ \hline
0.2        & Dev                     & 0.294396     & 0.294551 & 0.0524528       \\ \hline
0.1        & Dev                     & 0.353829     & 0.354145 & 0.0894466       \\ \hline
\end{tabular}
\caption{Exact and the winning DeVylder approximate values  of $\Phi_q$ for $f(x)=\frac{2}{3}e^{-x}+\frac{2}{3} e^{-2 x}$ when $\theta$ varies.}
\label{table:MixExp12Phiq}
\end{table}

\begin{table}[H]
\begin{tabular}{|l|l|l|l|l|}
\hline
$\theta$ & Closest approximation & Barrier exact & Barrier approx & \% error Barrier \\ \hline
1          & Dev                     & 3.45398         & 3.48756  & 0.972251         \\ \hline
0.9        & Dev                     & 3.20191         & 3.23103  & 0.909487         \\ \hline
0.8        & Dev                     & 2.90951         & 2.93074  & 0.729628         \\ \hline
0.7        & Dev                     & 2.57043         & 2.57742  & 0.272088         \\ \hline
0.6        & Dev                     & 2.1804          & 2.16054  & 0.910666         \\ \hline
0.5        & Ren                     & 1.74216         & 1.75266  & 0.60278          \\ \hline
0.4        & Expo                    & 1.2735          & 1.29456  & 1.65378          \\ \hline
0.3        & Expo                    & 0.81068         & 0.652264 & 19.5412          \\ \hline
\end{tabular}
\caption{Exact and approximate values of $b_{DeF}$ for for $f(x)=\frac{2}{3}e^{-x}+\frac{2}{3} e^{-2 x}$ when $\theta$ varies. As $b_{DeF}$ approaches $0$, errors of all the approximations increase dramatically, with the pure exponential approximation performing better than the rest. Meanwhile, as $b_{DeF}$ increases, Renyi and subsequently the \deV\ take the upper hand. For $\theta =0.2$ and $\theta =0.1$, all the approximations yield a 0 barrier approximate for exact barrier values of $0.392105$ and $0.0354538$ respectively, hence failing to predict the non-zero barrier.}
\label{table:MixExp12Bar}
\end{table}

\begin{table}[H]
\begin{tabular}{|l|l|l|l|l|l|}
\hline
    & \multicolumn{5}{c|}{J0}                                                                          \\ \hline
$\theta$    & J0 exact & J0   expo pure & J0   expo pure error & J0   expo CI & J0   expo CI error             \\ \hline
1   & 5.95034  & 5.99151        & 0.691856             & 6.26009      & 5.20551                        \\ \hline
0.9 & 5.15579  & 5.17573        & 0.386663             & 5.45269      & 5.7584                         \\ \hline
0.8 & 4.39383  & 4.38494        & 0.202205             & 4.67042      & 6.29489                        \\ \hline
0.7 & 3.68299  & 3.63933        & 1.18555              & 3.92937      & 6.68958                        \\ \hline
0.6 & 3.04577  & 2.96728        & 2.57704              & 3.25112      & 6.7423                         \\ \hline
0.5 & 2.50331  & 2.39942        & 4.15022              & 2.65901      & 6.21974                        \\ \hline
0.4 & 2.06833  & 1.9585         & 5.31006              & 2.17044      & 4.93725                        \\ \hline
0.3 & 1.74095  & 1.65616        & 4.86984              & 1.78984      & 2.80878                        \\ \hline
0.2 & 1.50439  & 1.44242        & 4.11969              & 1.50871      & 0.286672                       \\ \hline
0.1 & 1.30271  & 1.25324        & 3.79748              & 1.30271      & 0 \\ \hline
\end{tabular}
\caption{Values of $J_0$ compared with approximations using all exponential inputs ($J_0$ expo pure) and actual inputs but computed using the exponential formula ($J_0$ expo CI). The pure exponential approximation does a good job of approximating $J_0$ for higher values of $\theta$ considered, while the exponential CI approximation seemed to fair better for lower $\theta$ values}
\end{table}

\begin{table}[H]
\begin{tabular}{|l|l|l|l|l|l|}
\hline
    & \multicolumn{5}{c|}{a}                                                                        \\ \hline
$\theta$    & a exact  & a   expo pure & a   expo pure error & a   expo CI & a   expo CI error              \\ \hline
1   & 3.9669   & 3.99434       & 0.691861            & 4.17339     & 5.20551                        \\ \hline
0.9 & 3.4372   & 3.45049       & 0.386665            & 3.63512     & 5.7584                         \\ \hline
0.8 & 2.92922  & 2.9233        & 0.202204            & 3.11361     & 6.29489                        \\ \hline
0.7 & 2.45533  & 2.42622       & 1.18555             & 2.61958     & 6.68958                        \\ \hline
0.6 & 2.03051  & 1.97818       & 2.57704             & 2.16741     & 6.7423                         \\ \hline
0.5 & 1.66888  & 1.59961       & 4.15022             & 1.77268     & 6.21974                        \\ \hline
0.4 & 1.37888  & 1.30566       & 5.31006             & 1.44696     & 4.93725                        \\ \hline
0.3 & 1.16063  & 1.10411       & 4.86983             & 1.19323     & 2.80878                        \\ \hline
0.2 & 1.00293  & 0.961612      & 4.11969             & 1.0058      & 0.286672                       \\ \hline
0.1 & 0.868476 & 0.835496      & 3.79748             & 0.868476    &0 \\ \hline
\end{tabular}
\caption{Values of $a$ compared with approximations using all exponential inputs ($a$ expo pure) and actual inputs but computed using the exponential formula ($a$ expo CI)}
\end{table}

\begin{table}[H]
\begin{tabular}{|l|l|l|l|l|l|}
\hline
    & \multicolumn{5}{c|}{b}                                                                                                                                 \\ \hline
$\theta$    & b exact                        & b   expo pure                  & b   expo pure error & b   expo CI                    & b   expo CI error             \\ \hline
1   & 1.41036                        & 1.46188                        & 3.65293             & 1.25374                        & 11.1045                       \\ \hline
0.9 & 1.37645                        & 1.44439                        & 4.93621             & 1.23362                        & 10.3761                       \\ \hline
0.8 & 1.31492                        & 1.40417                        & 6.78809             & 1.19529                        & 9.09781                       \\ \hline
0.7 & 1.21057                        & 1.32258                        & 9.25207             & 1.12775                        & 6.84178                       \\ \hline
0.6 & 1.04634                        & 1.17215                        & 12.0245             & 1.01753                        & 2.7529                        \\ \hline
0.5 & 0.810767                       & 0.920406                       & 13.5229             & 0.853397                       & 5.25805                       \\ \hline
0.4 & 0.510085                       & 0.538725                       & 5.61475             & 0.634716                       & 24.4335                       \\ \hline
0.3 & 0.17425                        & 0.0105496                      & 93.9457             & 0.376872                       & 116.282                       \\ \hline
0.2 & 0 & 0 & 0           & 0.105322                       & 100 \\ \hline
0.1 & 0 & 0 & 0             & 0 & 0 \\ \hline
\end{tabular}
\caption{Values of $b$ compared with approximations using all exponential inputs ($b$ expo pure) and actual inputs but computed using the exponential formula ($b$ expo CI)}
\end{table}

\subsection{A Cram\'{e}r-Lundberg process with hyperexponential claims of order 3} \label{e:MixExp83}
Consider a Cram\'{e}r-Lundberg process with density function
$f(x)=\frac{12}{83 }e^{-x}+\frac{42}{83} e^{-2 x}+\frac{150}{83}e^{-3x}$, and $c=1$,  $\l=\frac{83}{48}$, $\th=\fr{263}{235}$, $p= \fr{263}{498}$, $q=\fr{5}{48}$.

The Laplace exponent of this process is
$\kappa(s) = s - \frac{12 s}{83 (s+1)}-\frac{21 s}{83 (s+2)}-\frac{50 s}{83 (s+3)}$ and from this one can invert $\frac{1}{\kappa(s) - q} =  \H{W}_q(s)$ to obtain the scale function \footnote{Laplace inversion done via Mathematica; coefficients and exponents are decimal approximations of the real values.}
\bea
W_q(x)  &= -0.0813294 e^{(-2.60997 x)} - 0.179472 e^{(-1.68854 x)} - 0.373887 e^{(-0.779311 x)}  + 1.63469 e^{(0.18198 x)}.
\eea
From this, we see that the dominant exponent is $\Phi_q = 0.18198$.

Figure \ref{fig:MixExp83} shows the exact and  approximate plots  of the first two derivatives of $W_q$. The exact plots    are labelled Wxexact, and coloured as the darkest.  The plots of $W_q'$  exhibit  noticeable unique minima around $x=2$, with the exact one being at $b_{DeF}=1.89732$, which  is the optimal barrier that  maximizes dividends here. Note that the  approximations are practically indistinguishable from the exact around this point (which is our main object of interest here).
\begin{figure}[H]
    \centering
    \begin{subfigure}[b]{0.45\textwidth}
        \includegraphics[width=\textwidth]{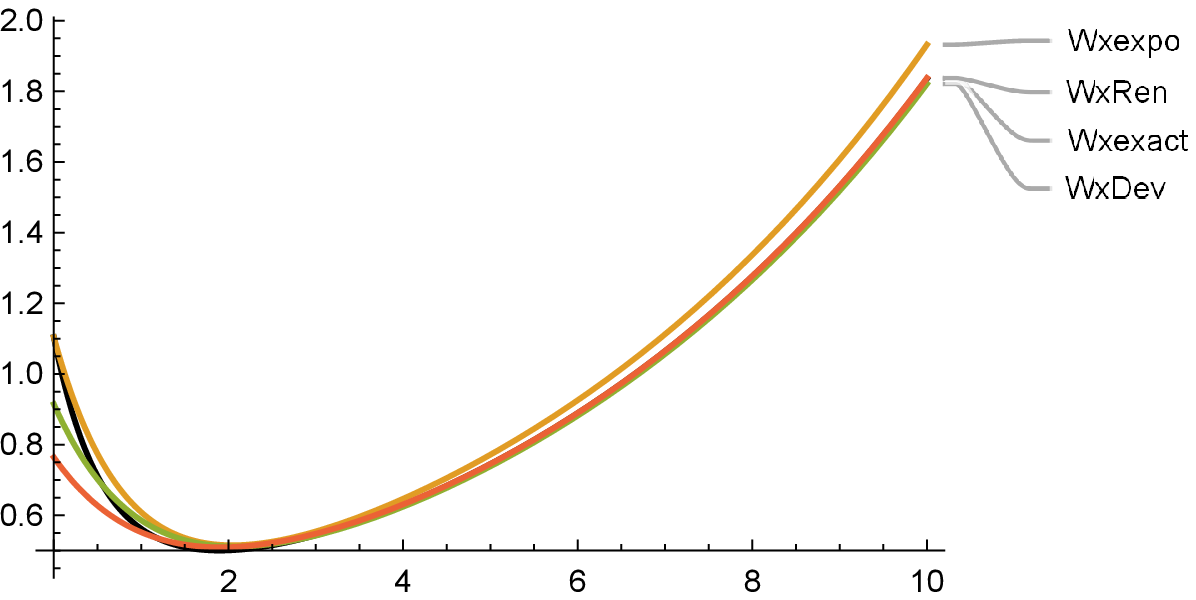}
        \caption{$W'_q(x)$}
        \label{fig:MixExp83W1}
    \end{subfigure}
    ~
    \begin{subfigure}[b]{0.45\textwidth}
        \includegraphics[width=\textwidth]{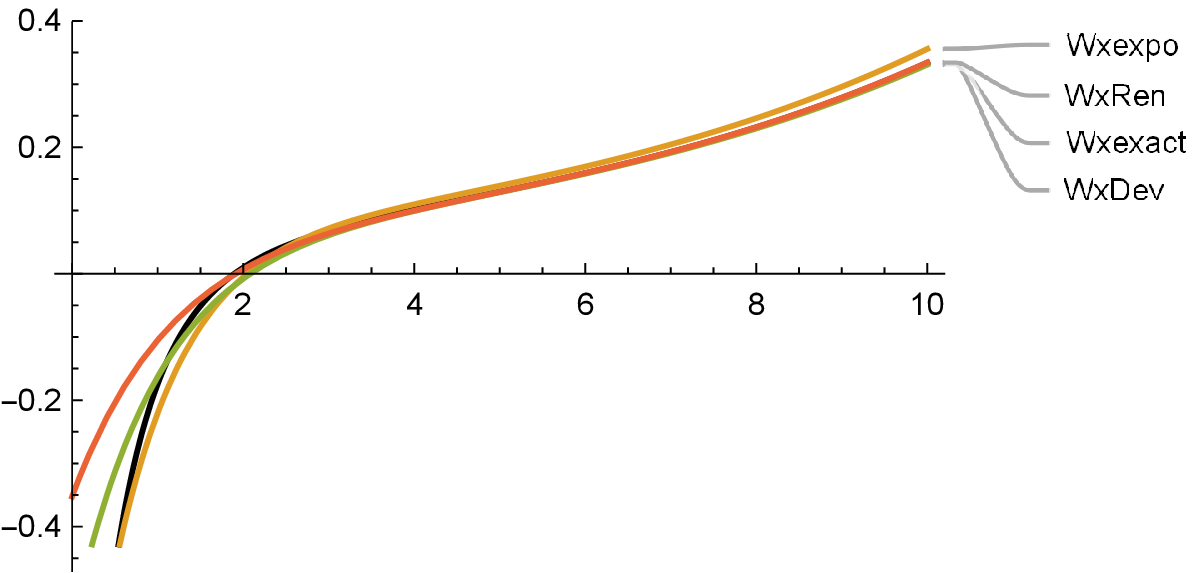}
        \caption{$W''_q(x)$}
        \label{fig:MixExp83W2}
    \end{subfigure}
    \caption{Exact  and  approximate plots of $W'_q(x)$ and $W''_q(x)$  for $f(x)=\frac{12}{83 }e^{-x}+\frac{42}{83} e^{-2 x}+\frac{150}{83}e^{-3x}$, $c=1$, $q=\fr{5}{48}$.}\label{fig:MixExp83}
\end{figure}
\begin{table}[H]
\begin{tabular}{|l|l|l|l|l|}
\hline
       & \begin{tabular}[c]{@{}l@{}}Dominant   exponent \\ $\Phi_q$\end{tabular} & \begin{tabular}[c]{@{}l@{}}Percent   relative error\\ ($\Phi_q$)\end{tabular} & \begin{tabular}[c]{@{}l@{}}Optimal barrier\\ $b_{DeF}$\end{tabular} & \begin{tabular}[c]{@{}l@{}}Percent   relative error\\ ($b_{DeF}$)\end{tabular} \\ \hline
Exact  & 0.18198                                                               & 0                                                                           & 1.89732                                                      & 0                                                                       \\ \hline
Expo   & 0.184095                                                              & 1.162215628                                                                 & 2.04608                                                      & 7.840532962                                                             \\ \hline
Renyi  & 0.181708                                                              & 0.149466974                                                                 & 2.08136                                                      & 9.699997892                                                             \\ \hline
Dev    & 0.182011                                                              & 0.017034839                                                                 & 1.91233                                                      & 0.79111589                                                               \\ \hline
\end{tabular}
\caption{Exact and approximate values of $\Phi_q$ and $b_{DeF}$ for $f(x)=\frac{12}{83 }e^{-x}+\frac{42}{83} e^{-2 x}+\frac{150}{83}e^{-3x}$, $c=1$, $q=\fr{5}{48}$, as well as percent relative errors, computed as the absolute value of the difference between the approximation and the exact, divided by the exact, times 100.
Relative errors for the $\Phi_q$ value are less than $1.25\%$, even for  the worse natural exponential  approximation, and the DeVylder approximation is the winner. The optimal barrier $b_{DeF}$ is  well approximated only by DeVylder.}
\label{table:MixExp83}
\end{table}

\begin{table}[H]
\begin{tabular}{|l|l|l|l|l|}
\hline
$\theta$ & Closest approximation & $\Phi_q$   exact & $\Phi_q$ approximation & \% error   $\Phi_q$ \\ \hline
263/235    & Dev                   & 0.18198        & 0.182011      & 0.0168217         \\ \hline
243/235    & Dev                   & 0.194712       & 0.194754      & 0.0213671         \\ \hline
223/235    & Dev                   & 0.209221       & 0.209279      & 0.0274827         \\ \hline
203/235    & Dev                   & 0.225876       & 0.225957      & 0.0358309         \\ \hline
183/235    & Dev                   & 0.245146       & 0.245262      & 0.0474032         \\ \hline
163/235    & Dev                   & 0.267635       & 0.267806      & 0.0637063         \\ \hline
143/235    & Dev                   & 0.294126       & 0.294382      & 0.0870647         \\ \hline
123/235    & Dev                   & 0.325643       & 0.326038      & 0.121115          \\ \hline
103/235    & Dev                   & 0.363539       & 0.364163      & 0.171618          \\ \hline
83/235     & Dev                   & 0.40961        & 0.410625      & 0.247788          \\ \hline
63/235     & Dev                   & 0.466261       & 0.46796       & 0.364457          \\ \hline
43/235     & Dev                   & 0.536719       & 0.539647      & 0.545532          \\ \hline
23/235     & Dev                   & 0.62533        & 0.630516      & 0.829419          \\ \hline
3/235      & Dev                   & 0.737962       & 0.747389      & 1.27736           \\ \hline
\end{tabular}
\caption{Exact and the winning DeVylder  approximate values  of $\Phi_q$, for $f(x)=\frac{12}{83 }e^{-x}+\frac{42}{83} e^{-2 x}+\frac{150}{83}e^{-3x}$, $\theta$ varies.}
\label{table:MixExp83Phiq}
\end{table}

\begin{table}[H]
\begin{tabular}{|l|l|l|l|l|}
\hline
$\theta$ & Closest approximation & Barrier exact & Barrier approx & \% error Barrier \\ \hline
263/235 & Dev    & 1.89732   & 1.91233  & 0.791183 \\ \hline
243/235 & Dev    & 1.79954   & 1.78002  & 1.08482  \\ \hline
183/235 & Ren    & 1.45224   & 1.52484  & 4.9989   \\ \hline
163/235 & Ren    & 1.31579   & 1.33691  & 1.60463  \\ \hline
143/235 & Ren    & 1.16804   & 1.12368  & 3.79796  \\ \hline
123/235 & Expo   & 1.00898   & 1.04123  & 3.19653  \\ \hline
103/235 & Expo   & 0.839228  & 0.794964 & 5.27444  \\ \hline
83/235  & Expo   & 0.660338  & 0.513179 & 22.2854  \\ \hline
63/235  & Expo   & 0.474896  & 0.196234 & 58.6785  \\ \hline
43/235  & Expo   & 0.286563  & 0        & 100      \\ \hline
23/235  & Expo   & 0.0998863 & 0        & 100      \\ \hline
3/235   & All  & 0         & 0        & 0        \\ \hline
\end{tabular}
\caption{Exact and approximate values of $b_{DeF}$ for $f(x)=\frac{12}{83 }e^{-x}+\frac{42}{83} e^{-2 x}+\frac{150}{83}e^{-3x}$, when $\theta$ varies. Unsurprisingly, when $b_{DeF}$ is near $0$, the natural exponential approximation  wins, but has poor performance, and as $b_{DeF}$ increases, Renyi and subsequently the \deV\ take the upper hand. For the smallest three values of $\theta$, all of the approximations yielded $\T b_{DeF}=0$ as the optimal barrier, with this  being true only for $\theta=3/235$.}
\label{table:MixExp83Bar}
\end{table}

We move now to the dividend problem with capital injections with cost $k \geq 1$ as in Theorem \ref{p:cost}. One can compute  the value function  $J_0$ at $x=0$ in terms of $W$, $Z$,$C$, $S$, and $G$ -- see equation \ref{struct}.

To provide a more concrete example, fixing $q=\fr{5}{48}$,$P=0$, $k=3/2$ as input parameters we compute for values of $J_0$ as a function of $\theta$, with results summarized in the tables below. The tables provide comparisons of the computed optimal quantities $J_0$, $a$, and $b$ to an approximation using all exponential inputs (referred to as $J_0$, $a$, and $b$ expo pure) and to an approximation which uses actual inputs but computed using the exponential formula as described in equation \ref{costexp} (referred to as $J_0$, $a$, and $b$ expo CI).

\begin{table}[H]
\begin{tabular}{|l|l|l|l|l|l|}
\hline
           & \multicolumn{5}{c|}{J0}                                                                            \\ \hline
$\theta$ & $J_0$   & $J_0$ expo pure & $J_0$ expo pure error & $J_0$ expo CI & $J_0$ expo CI error             \\ \hline
263/235    & 3.7747     & 3.76883        & 0.155556             & 4.11784      & 9.09041                        \\ \hline
243/235    & 3.41491    & 3.38603        & 0.845606             & 3.74156      & 9.5654                         \\ \hline
223/235    & 3.0636     & 3.00802        & 1.81444              & 3.36985      & 9.99637                        \\ \hline
203/235    & 2.72335    & 2.63828        & 3.1238               & 3.00466      & 10.3296                        \\ \hline
183/235    & 2.39737    & 2.28225        & 4.80185              & 2.64879      & 10.4871                        \\ \hline
163/235    & 2.08958    & 1.94765        & 6.79226              & 2.3062       & 10.3665                        \\ \hline
143/235    & 1.80446    & 1.64396        & 8.8946               & 1.9823       & 9.85516                        \\ \hline
123/235    & 1.54668    & 1.38072        & 10.73                & 1.68379      & 8.86472                        \\ \hline
103/235    & 1.32041    & 1.16526        & 11.7499              & 1.4178       & 7.37587                        \\ \hline
83/235     & 1.12864    & 1.00194        & 11.2258              & 1.19022      & 5.45555                        \\ \hline
63/235     & 0.972835   & 0.88785        & 8.73585              & 1.00404      & 3.20798                        \\ \hline
43/235     & 0.852739   & 0.789923       & 7.36635              & 0.859039     & 0.738837                       \\ \hline
23/235     & 0.751597   & 0.701299       & 6.69218              & 0.751597     & 0 \\ \hline
3/235      & 0.660372   & 0.620567       & 6.02761              & 0.660372     & 0 \\ \hline
\end{tabular}
\caption{Values of $J_0$ compared with approximations using all exponential inputs ($J_0$ expo pure) and actual inputs but computed using the exponential formula ($J_0$ expo CI). The pure exponential approximation does a good job of approximating $J_0$ for higher values of $\theta$ considered, while the exponential CI approximation seemed to fair better for lower $\theta$ values}
\end{table}

\begin{table}[H]
\begin{tabular}{|l|l|l|l|l|l|}
\hline
           & \multicolumn{5}{c|}{a}                                                                        \\ \hline
$\theta$ & $a$  & $a$   expo pure & $a$   expo pure error & $a$   expo CI & $a$   expo CI error              \\ \hline
263/235    & 2.51647  & 2.74523       & 0.155536            & 2.51255     & 9.09042                        \\ \hline
243/235    & 2.27661  & 2.49437       & 0.845596            & 2.25735     & 9.56541                        \\ \hline
223/235    & 2.0424   & 2.24657       & 1.81443             & 2.00535     & 9.99638                        \\ \hline
203/235    & 1.81557  & 2.00311       & 3.1238              & 1.75885     & 10.3296                        \\ \hline
183/235    & 1.59825  & 1.76586       & 4.80185             & 1.5215      & 10.4871                        \\ \hline
163/235    & 1.39306  & 1.53747       & 6.79226             & 1.29844     & 10.3665                        \\ \hline
143/235    & 1.20298  & 1.32153       & 8.8946              & 1.09598     & 9.85516                        \\ \hline
123/235    & 1.03112  & 1.12252       & 10.73               & 0.920479    & 8.86472                        \\ \hline
103/235    & 0.880271 & 0.945198      & 11.7499             & 0.77684     & 7.37587                        \\ \hline
83/235     & 0.752428 & 0.793477      & 11.2258             & 0.667962    & 5.45555                        \\ \hline
63/235     & 0.648557 & 0.669362      & 8.73585             & 0.5919      & 3.20798                        \\ \hline
43/235     & 0.568493 & 0.572693      & 7.36635             & 0.526616    & 0.738838                       \\ \hline
23/235     & 0.501065 & 0.501065      & 6.69218             & 0.467532    & 0 \\ \hline
3/235      & 0.440248 & 0.440248      & 6.02761             & 0.413711    & 0  \\ \hline
\end{tabular}
\caption{Values of $a$ compared with approximations using all exponential inputs ($a$ expo pure) and actual inputs but computed using the exponential formula ($a$ expo CI)}
\end{table}

\begin{table}[H]
\begin{tabular}{|l|l|l|l|l|l|}
\hline
           & \multicolumn{5}{c|}{b}                                                                                                                                 \\ \hline
$\theta$ & $b$                        & $b$   expo pure                  & $b$   expo pure error & $b$   expo CI                    & $b$   expo CI error             \\ \hline
263/235    & 0.709355                       & 0.805116                       & 13.4997             & 0.677918                       & 4.43179                       \\ \hline
243/235    & 0.695874                       & 0.801936                       & 15.2416             & 0.671779                       & 3.46259                       \\ \hline
223/235    & 0.677601                       & 0.794377                       & 17.2337             & 0.662801                       & 2.18425                       \\ \hline
203/235    & 0.653005                       & 0.779265                       & 19.3352             & 0.649805                       & 0.490147                      \\ \hline
183/235    & 0.620126                       & 0.751601                       & 21.2012             & 0.631097                       & 1.76912                       \\ \hline
163/235    & 0.576553                       & 0.704104                       & 22.123              & 0.604293                       & 4.81128                       \\ \hline
143/235    & 0.519526                       & 0.627369                       & 20.7579             & 0.566198                       & 8.98366                       \\ \hline
123/235    & 0.446259                       & 0.511076                       & 14.5246             & 0.512961                       & 14.947                        \\ \hline
103/235    & 0.354524                       & 0.346046                       & 2.39143             & 0.440755                       & 24.3231                       \\ \hline
83/235     & 0.243362                       & 0.126054                       & 48.2032             & 0.347059                       & 42.6099                       \\ \hline
63/235     & 0.113593                       & 0 & 100              & 0.231975                       & 104.216                       \\ \hline
43/235     & 0 & 0 & 0             & 0.0987484                      & 0 \\ \hline
23/235     & 0 & 0 & 0             & 0 & 0                       \\ \hline
3/235      & 0 & 0 & 0             & 0 & 0 \\ \hline
\end{tabular}
\caption{Values of $b$ compared with approximations using all exponential inputs ($b$ expo pure) and actual inputs but computed using the exponential formula ($b$ expo CI)}
\end{table}

To provide a  point of comparison, we fix $q=\fr{5}{48}$, and compute the de Finetti barrier to be $b_{DeF} = 1.89732$ and the corresponding dividend value function when starting at $x=0$ to be $J_{DeF} = 1.99847$.

\begin{table}[H]
\begin{tabular}{|l|l|l|l|l|l|}
\hline
$k$     & $J_0$ \% deviation & $J_0 - J_{DeF}$ & $J_0$ & $b$ \% deviation & $b$                         \\ \hline
1     & 154.123      & 3.0801          & 5.07857 & 100          & 0 \\ \hline
2     & 65.714       & 1.31327         & 3.31174 & 43.0208      & 1.08108                        \\ \hline
3     & 42.863       & 0.856604        & 2.85507 & 25.4792      & 1.4139                         \\ \hline
4     & 31.4465      & 0.628448        & 2.62692 & 17.6853      & 1.56178                        \\ \hline
5     & 24.6995      & 0.493611        & 2.49208 & 13.4234      & 1.64264                        \\ \hline
6     & 20.2855      & 0.4054          & 2.40387 & 10.7759      & 1.69287                        \\ \hline
7     & 17.1884      & 0.343504        & 2.34197 & 8.98441      & 1.72686                        \\ \hline
8     & 14.9014      & 0.2978          & 2.29627 & 7.69619      & 1.7513                         \\ \hline
9     & 13.1463      & 0.262724        & 2.26119 & 6.72732      & 1.76968                        \\ \hline
10    & 11.758       & 0.23498         & 2.23345 & 5.97306      & 1.78399                        \\ \hline
100   & 1.11095      & 0.0222019       & 2.02067 & 0.533448     & 1.8872                         \\ \hline
1000  & 0.110409     & 0.00220648      & 2.00067 & 0.0527257    & 1.89632                        \\ \hline
10000 & 0.0109536    & 0.000218903     & 1.99869 & 0.00533526   & 1.89722                        \\ \hline
\end{tabular}
\caption{Values of $J_0$ and $b$ in presence of capital injections compared to the case where capital injections are non-existent, $J_{DeF} = 1.99847$ and $b_{DeF} = 1.89732$. As $k$ is increased one can see that $J_0$ and $b$ approaches $J_{DeF}$ and $b_{DeF}$. This is expected since higher costs of injecting capital makes it less viable, hence it is treated like the concept does not exist.}
\label{table:MixExp83JDef}
\end{table} 
\subsection{A Cram\'{e}r-Lundberg process with oscillating density and scale function} \label{e:MatExpCos}

In the following example, we study a \CL\ model with density of claims given by
\begin{align*}
f(x)=& u e^{-a x} 2 \cos ^{2}\left(\frac{\omega x+\phi}{2}\right)=u e^{-a x}(1+\cos (\omega x+\phi))=\\=& e^{-a x}(u+u \cos (\phi) \cos (\omega x)-u \sin (\phi) \sin (\omega x))
\end{align*}
where
\[u=\frac{a\left(a^{2}+\omega^{2}\right)}{a^{2}+\omega^{2}+a^{2} \cos (\phi)-a \omega \sin (\phi)}. \]
%One can check that $\int f(x) d x=1$ with such value of $u$.

Assuming further that $a=1$, $\phi=2$, $\omega=20$, and that $\theta=1$, $q=1/10$, the Laplace exponent for this process is
$\kappa(s) = \frac{s \left(2.09898 s^3+5.29695 s^2+843.502 s+420.846\right)}{(s+1.) \left(s^2+2. s+401.\right)}$ and the scale function is
\begin{align*}
W_q(x)  &= 0.824723 e^{0.0881484 x} -0.348141 e^{-0.540677 x} \\
   &+e^{-1.0117 x} \cos (19.9957 x) ~\Big( -(0.000285494\, +0.0000804151 i) \sin (39.9914 x)\\
   &-(0.0000804151\, +0.000285494 i)+(-0.0000804151+0.000285494 i) \cos (39.9914 x) \Big)\\
   &+e^{-1.0117 x} \sin (19.9957 x) ~\Big( -(0.0000804151\, -0.000285494 i) \sin (39.9914 x) \\
   & +(0.000285494\, +0.0000804151 i) \cos (39.9914 x)-(0.000285494\, -0.0000804151 i) \Big).
\end{align*}

\begin{figure}[H]
    \centering
    \begin{subfigure}[b]{0.45\textwidth}
        \includegraphics[width=\textwidth]{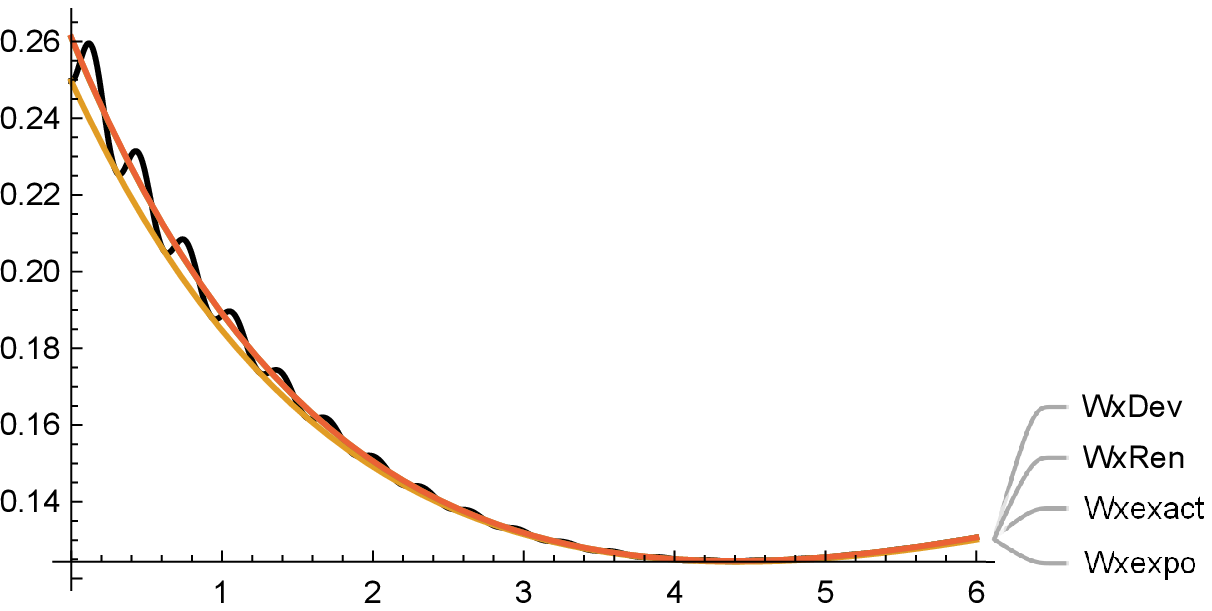}
        \caption{$W'_q(x)$}
        \label{fig:MatExpCosW1}
    \end{subfigure}
    ~
    \begin{subfigure}[b]{0.45\textwidth}
        \includegraphics[width=\textwidth]{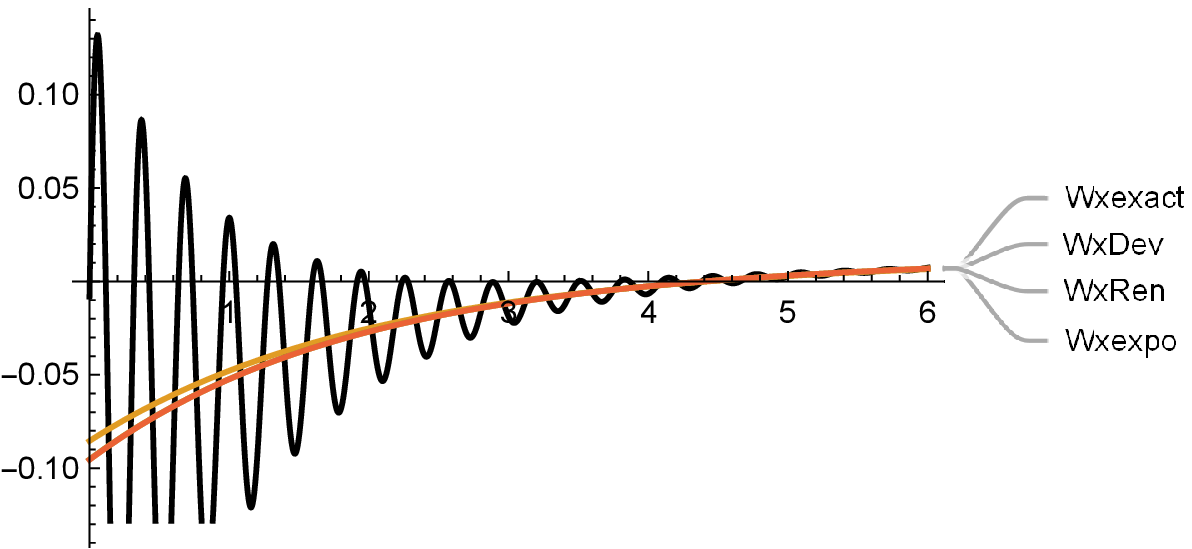}
        \caption{$W''_q(x)$}
        \label{fig:MatExpCosW2}
    \end{subfigure}
    \caption{Plots of $W'_q(x)$, and $W''_q(x)$ of the exact solution and the approximations for $f(x)= u e^{-a x} 2 \cos ^{2}\left(\frac{\omega x+\phi}{2}\right)$, $\th =1$, $q=\fr{1}{10}$.}\label{fig:MatExpCosW}
\end{figure}

\begin{table}[!h]
\begin{tabular}{|l|l|l|l|l|}
\hline
       & \begin{tabular}[c]{@{}l@{}}Dominant   exponent \\ $\Phi_q$\end{tabular} & \begin{tabular}[c]{@{}l@{}}Percent   relative error\\ ($\Phi_q$)\end{tabular} & \begin{tabular}[c]{@{}l@{}}Optimal barrier\\ $b_{DeF}$\end{tabular} & \begin{tabular}[c]{@{}l@{}}Percent   relative error\\ ($b_{DeF}$)\end{tabular} \\ \hline
Exact  & 0.0881484                  & 0                                 & 4.38201              & 0                             \\ \hline
Expo   & 0.0878658                  & 0.32053                           & 4.42263              & 0.927122                      \\ \hline
Renyi  & 0.0881481                  & 0.000314617                       & 4.39788              & 0.362284                      \\ \hline
Dev    & 0.0881484                  & 6.11743*10\textasciicircum{}-6    & 4.39745              & 0.352331                      \\ \hline
\end{tabular}
\caption{Exact and approximate values of $\Phi_q$ and $b_{DeF}$ for $f(x)= u e^{-a x} 2 \cos ^{2}\left(\frac{\omega x+\phi}{2}\right)$, $\th =1$, $q=\fr{1}{10}$. The DeVylder approximation wins on both fronts.}
\label{table:MatExpCos}
\end{table}

Clearly, 
 our completely monotone approximation cannot fully reproduce functions like $W_q'(x), W_q''(x)$ in examples like this where  oscillations occur (note however that the \deF\ optimal barrier is well approximated here). If a more exact reproduction is  necessary,
 higher order approximations should be used.

%\newpage 

%\input{exa} \input{315}
%\input{NH5} %bad graphs
%\input{NH5bic} %bad graphs
%\input{exaRein}
%\input{Exp235per}
\section{The maximal error of exponential  \app s $J_0$ along  one parameter families of Cram\'{e}r-Lundberg processes} \label{e:Exp32per}

In this section, we provide the two  approximations for the dividend value with capital injections $J_0$, and the dividend barrier $b$, for two  one parameter families of Cram\'{e}r-Lundberg processes, with   densities  given \resp\ by:

\beq && \label{f3}
f(x)=k_\e \le[e^{-x} + \epsilon e^{-2 x}\ri]
 \\ &&
 \label{f2}
f(x)=k_\e \le[\frac{12}{83 }e^{-x}+ \epsilon \le(\frac{42}{83} e^{-2 x}+\frac{150}{83}e^{-3x}\ri)\ri],
\eeq
where $k_\e$ is the normalization constant, and compute the {maximal error of approximation} when $\e \in (0,\I)$ and $ \th \approx 1$.  For this choice, the pure exponential approximation works considerably better, {$\for \e$}.

\begin{table}[H]
\begin{tabular}{|l|l|l|l|l|l|}
\hline
      & \multicolumn{5}{c|}{J0}                                                              \\ \hline
      & J0 exact & J0   expo pure & J0   expo pure error & J0   expo CI & J0   expo CI error \\ \hline
0.001 & 7.1879   & 7.18802        & 0.0016603            & 7.18849      & 0.00827967         \\ \hline
0.01  & 7.17075  & 7.17193        & 0.0164663            & 7.17666      & 0.0824782          \\ \hline
0.1   & 7.008    & 7.01863        & 0.151653             & 7.06358      & 0.793041           \\ \hline
1     & 5.95034  & 5.99151        & 0.691856             & 6.26009      & 5.20551            \\ \hline
10    & 4.20175  & 4.19406        & 0.183122             & 4.40089      & 4.73941            \\ \hline
100   & 3.66909  & 3.6654         & 0.100631             & 3.69555      & 0.721228           \\ \hline
1000  & 3.6025   & 3.60208        & 0.0117065            & 3.60523      & 0.075585           \\ \hline
\end{tabular}
\caption{$\l=1$, $\th=1$, $q=\frac{1}{10}$ $k=3/2$ and $P=0$. As expected, the errors decrease both as $\epsilon$ goes to zero and infinity since the   densities approach an exponential density.}
\end{table}

\begin{figure}[H]
    \centering
    \begin{subfigure}[b]{0.45\textwidth}
        \includegraphics[width=\textwidth]{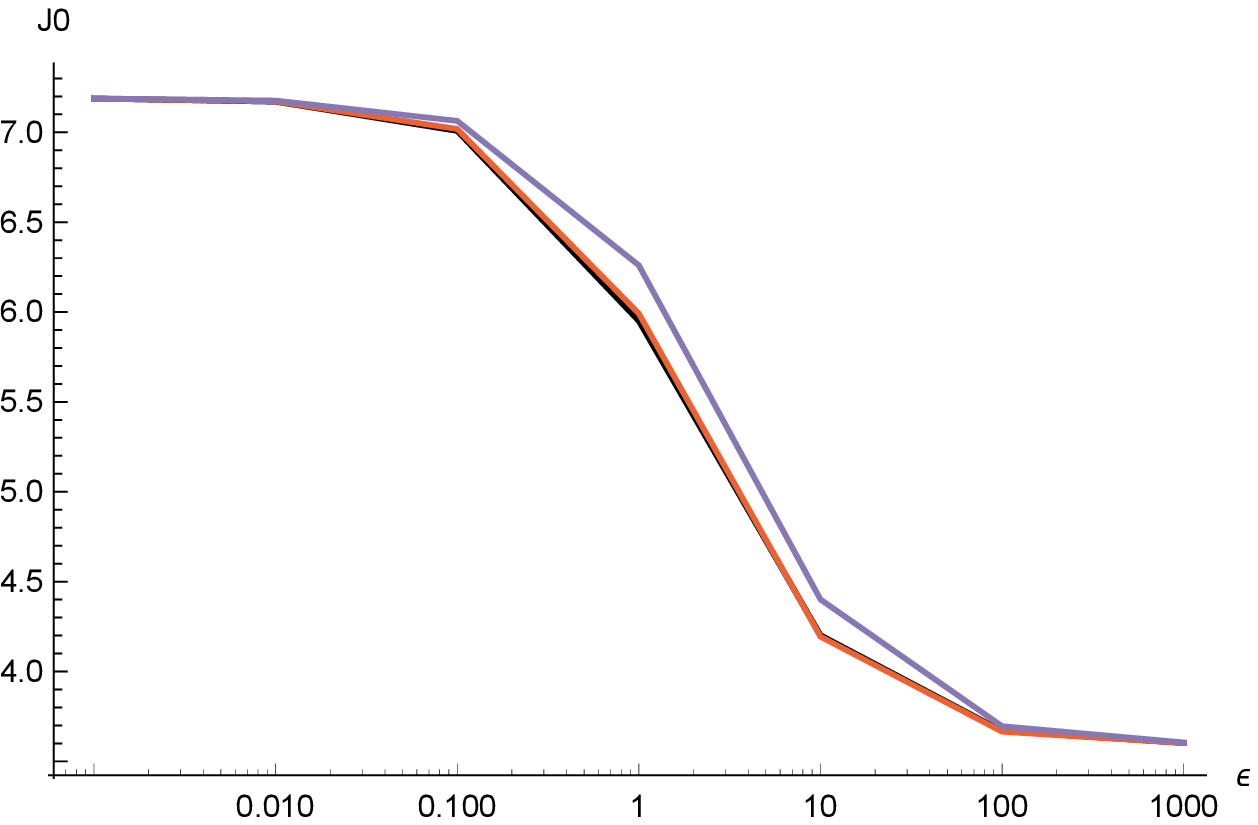}
        \label{fig:MixExp12J0eps}
    \end{subfigure}
    ~
    \begin{subfigure}[b]{0.45\textwidth}
        \includegraphics[width=\textwidth]{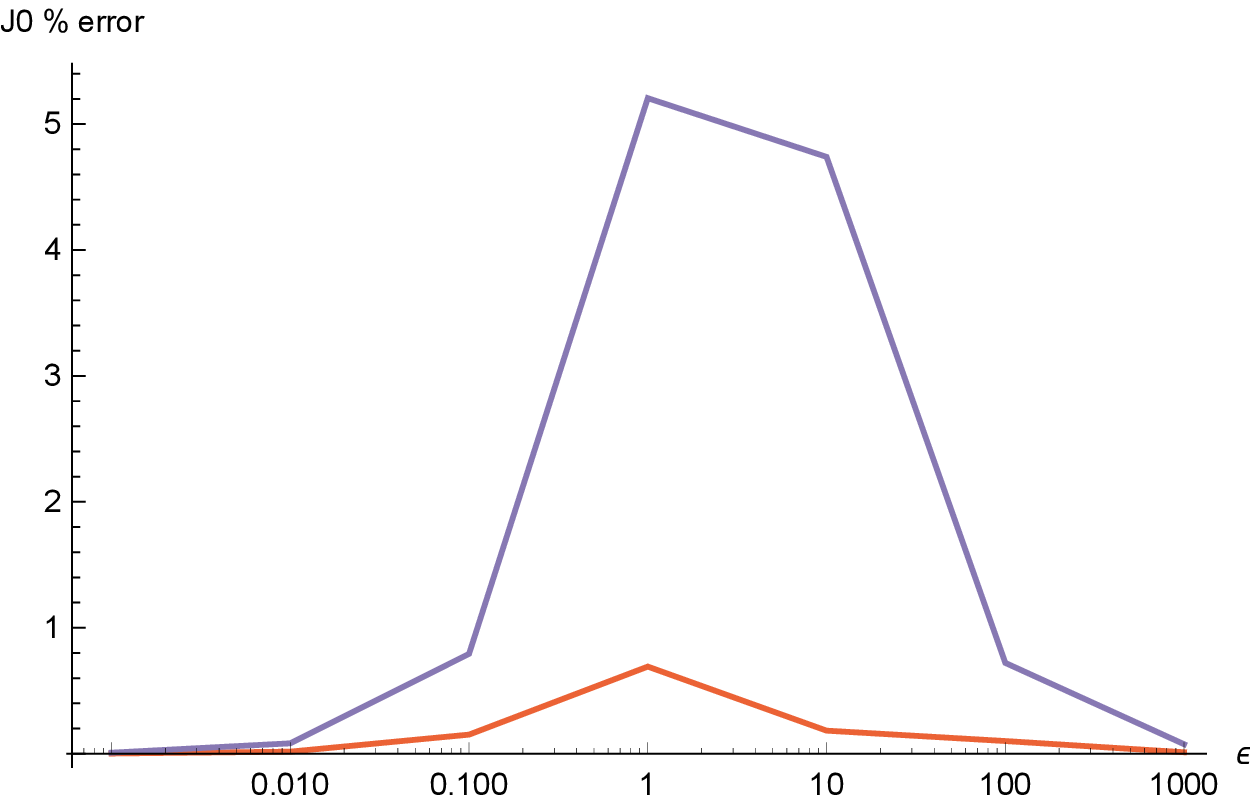}
        \label{fig:MixExp12J0epserr}
    \end{subfigure}
    \caption{$J_0$ values and errors plotted against $\epsilon$. Errors peak at $\epsilon=1$.}\label{fig:MixExp12J0}
\end{figure}

%\subsection{A family of densities of order three with $\th$ fixed}

We do the same thing  for the family of   densities given by $f(x)=k_\e \le[\frac{12}{83 }e^{-x}+ \epsilon \le(\frac{42}{83} e^{-2 x}+\frac{150}{83}e^{-3x}\ri)\ri]$.

\begin{table}[H]
\begin{tabular}{|l|l|l|l|l|l|}
\hline
      & \multicolumn{5}{c|}{J0}                                                              \\ \hline
$\epsilon$      & J0 exact & J0   expo pure & J0   expo pure error & J0   expo CI & J0   expo CI error \\ \hline
0.001 & 7.95508  & 7.95771        & 0.0330111            & 7.96565      & 0.132796           \\ \hline
0.01  & 7.68765  & 7.71127        & 0.307292             & 7.78772      & 1.30166            \\ \hline
0.1   & 6.06176  & 6.15381        & 1.5186               & 6.62641      & 9.31507            \\ \hline
1     & 3.7747   & 3.76883        & 0.155556             & 4.11784      & 9.09041            \\ \hline
10    & 3.1382   & 3.1379         & 0.00959692           & 3.23354      & 3.03813            \\ \hline
100   & 3.05894  & 3.06427        & 0.174306             & 3.12284      & 2.08921            \\ \hline
1000  & 3.0508   & 3.05678        & 0.196219             & 3.11149      & 1.98956            \\ \hline
\end{tabular}
\caption{$\l=1$, $c=1$, $q=\frac{5}{48}$ $k=3/2$ and $P=0$. As $\epsilon$ goes to zero, the density becomes exponential hence the decrease in errors. As $\epsilon$ goes to infinity, the density approaches a hyper exponential density of order 2, but still both methods of approximating $J_0$ yield reasonable results.}
\end{table}

\begin{figure}[H]
    \centering
    \begin{subfigure}[b]{0.45\textwidth}
        \includegraphics[width=\textwidth]{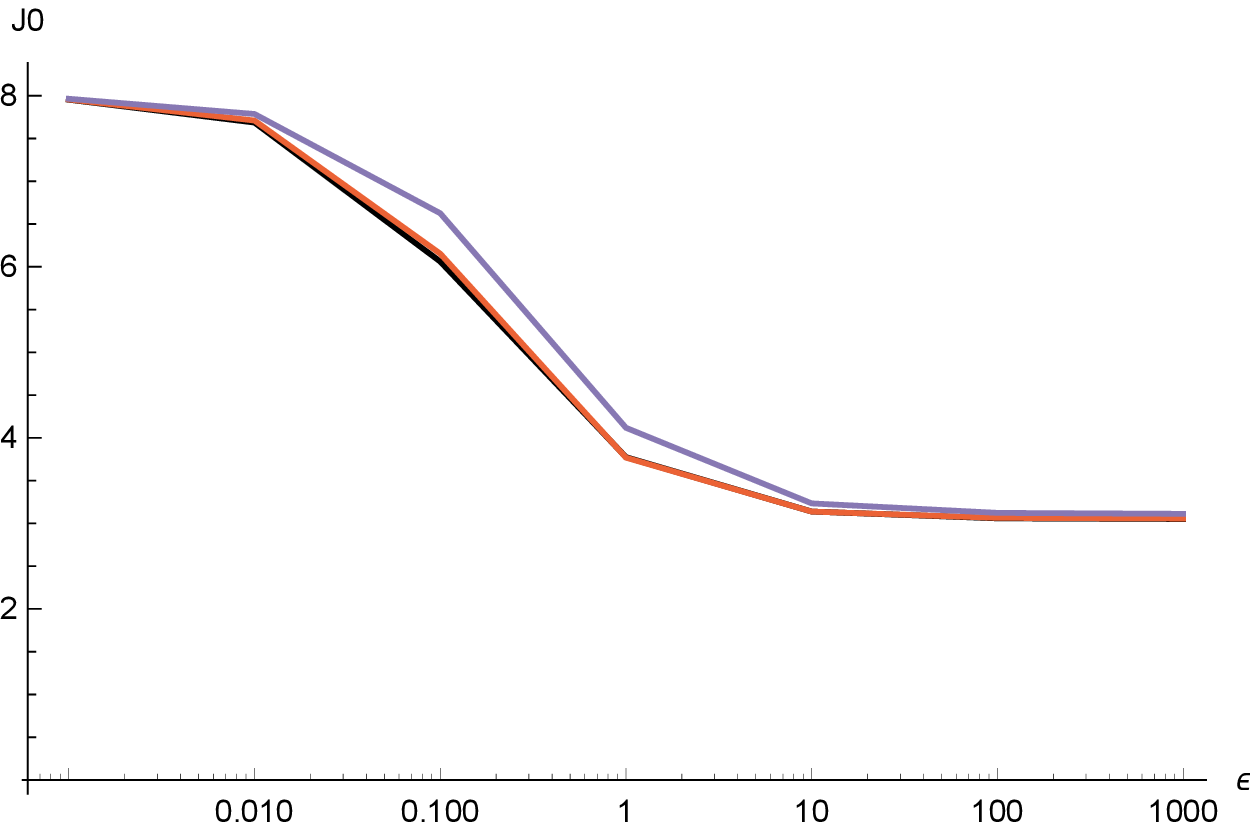}
        \label{fig:MixExp83J0eps}
    \end{subfigure}
    ~
    \begin{subfigure}[b]{0.45\textwidth}
        \includegraphics[width=\textwidth]{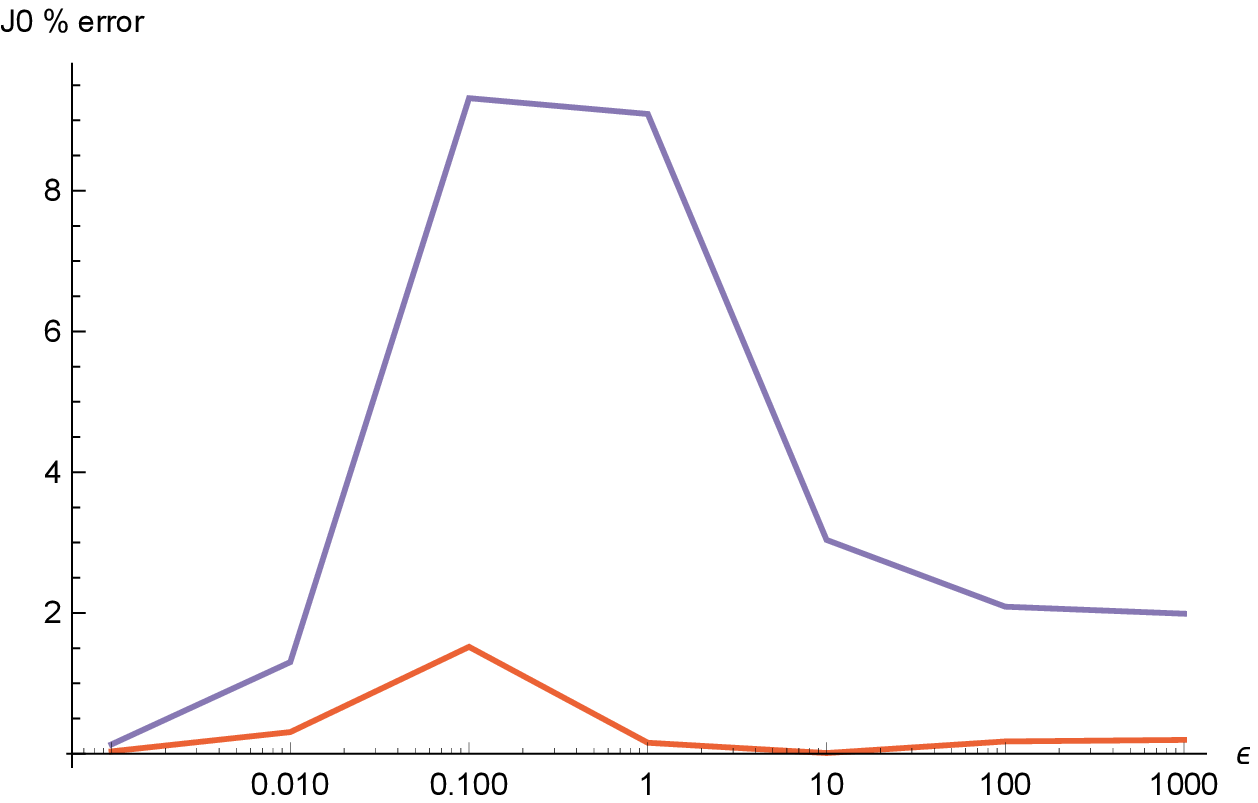}
        \label{fig:MixExp83J0epserr}
    \end{subfigure}
    \caption{$J_0$ values and errors plotted against $\epsilon$. Errors peak at $\epsilon=0.1$.}\label{fig:MixExp83J0}
\end{figure}

\sec{The profit function when  the claims are distributed according to a \me\ density \la{s:me}}

Consider now the more general case when the claims are distributed according to a matrix exponential density generated by a row vector $\vb$ and by an invertible matrix $B$ of order $n$, which are such that the vector $\vb e^{ x B}$ is decreasing componentwise to $0$, and $\vb . \vo \neq 0 $, with $\vo %, \bff b
$ a  column vector. As customary, we  restrict \wlo to the case
when  $\vb$ is a probability vector, and  $\vb . \bff 1 =1$, so that  $$\ovl F(x)=%(\vb . \vo)^{-1}
\vb e^{ x B} \vo$$  is a valid survival function.

 The matrix versions of our  functions are:

  \be \la{mpf} \bc C_a(x) =\l \int_0^x W_q(x-y) \;   \ovl F(y+a)  \; dy = \l \vb \int_0^x W_q(x-y) \;e^{ y B}     \; dy  \; e^{ a B} \vo= \vec C(x) e^{ a B} \vo
 \\ m_a(y)=\int_0^a z f(y+z)dz=\vb  \; e^{ y B} \; \int_0^a z \;e^{ z B} (-B)    \; dz \;\vo=\vb   \; e^{ y B} M(a) \vo
 \\G_a(x)=\l \int_0^x  W_q(x-y) \;  m_a(y) \; \; dy= \vec C(x) M(a) \vo \; \; \; \; \ec, \ee
 where
\be \bc C(x)=\l  \int_0^x W_q(x-y) \;e^{ y B}     \; dy\\
\vec C(x)=\l \vb \int_0^x W_q(x-y) \;e^{ y B}     \; dy \ec. \ee

The product  formulas \eqr{mpf} may also be established directly in the \PH case,  using the conditional independence of the ruin probability of the overshoot size. % -see section \ref{\eeXa

We derive first these extensions from scratch for $(\vb, B)$ \PH densities,
in order to highlight their probabilistic interpretation. Later, we will  show  that the \me\ case follows as a  particular case of \cite{Gaj}.

Recall first \cite{AA} that $\Rui_q(x)=\vRui_q(x) \vo,$
where  $\vRui_q(x)$ is  a vector whose components represent the probability that ruin occurs during a certain phase,
and that the conditional independence of ruin and overshoots translates into the product  formula
 \be \Rui_q(x,y):=P_x[\tz <\I, X_{\tz} < -y]= \vRui_q(x) e^{ y B} \vo.\ee

To take advantage of this, it is convenient to replace
from the beginning $Z_q(x)$ by $\Rui_q(x)$, taking advantage of the formula \cite{AKP,Kyp}
\be  Z_q(x)= \Rui_q(x)+  W_q(x) \frq \Lra C(x)=(c-\frq) W_q(x)- \Rui_q(x).\ee
%(see  -- basically, $ Z_q(x)$ is a ``version" of $\Rui_q(x)$ which fits smoothly the boundary condition $\Rui_q(x)=1, x<0$).

Alternatively, one may introduce a  vector function

 \be \vec Z_q(x):= \vec \Rui_q(x)+  W_q(x)  \frq \vo. \ee

\Oth the mean function \mbw $$m_a= \int_0^a y \; F( \md y)=- a \ovl F(a) + \int_0^a \ovl F(x) dx=\vb M(a) \vo, M(a)=-B^{-1} - e^{a B}\pr{a \C_n-B^{-1}}.$$

 The following result follows in the \PH case just as in the  exponential case \cite{AGLW}; in the \me\ case, it may be obtained from \cite{Gaj}:
\beP
For a % {\per}
\CL process  (\cP) with \me\ of type $(\vb,B)$, \ith
\BEN \im
 \be J_x= \bc k G_a(x) +J_0 S_a(x)=k G_a(x) +\fr{1-k G'(b)}{S'(b)} S_a(x),  &x \in[0, b]\\k x+J_0
&x \in[-a, 0]\\ 0 & x \leq -a \ec, \la{structPH}\ee
where
\be \bc C(x)=\l  \int_0^x W_q(x-y) \;e^{ y B}     \; dy\\
\vec C(x)=\l \vb \int_0^x W_q(x-y) \;e^{ y B}     \; dy=c W_q(x)\vec 1-\vec Z_q(x)=(c -\frq) W_q(x)\vec 1-\vec \Rui_q(x)\\  G_a(x)= \vec C_\q(x) M(a)\bff 1\\
R_a(x)=S_a(x) -  Z_q(x)=  \vec C_\q(x) e^{a B} \bff 1
 %{+k \si W_q(x), \text{ cf sols 2,3?}}
\ec, \ee
and
  \be J_0=
\fr{1-k \vec C_\q'(b)\; M(a) \bff 1}{q W_q(b) + \vec C_\q'(b) e^{a B} \vo}. \la{J0PH}
 \ee

 \im  For fixed $a$, the optimality equation $\fr{\partial}{\partial b} J_0^{a,b}
=0$ simplifies to
  \be J_0=
\fr{k \vec C_\q''(b)\; M(a) \bff 1}{q W_q'(b) + \vec C_\q''(b) e^{a B} \vo}.   \la{J0aPH}\ee
\EEN

\eeP

\iffalse
\prf
%These results follow from those of \cite{Gaj}.
1. It remains only to show the identity of the two formulas for $C(x)$. This holds \red{since ...}
\qed
\fi

\beR The additive separation of $a,b$ which was the basis of proving optimality in the exponential case does not  seem possible anymore, but
\eqr{J0PH} allows the numeric computation of  the optimum. \eeR
\iffalse 
\beR \red{To be removed eventually}. Note that  equations related to  \eqr{J0a},
\eqr{strEqa} are given in \cite[Cor. 12]{Gaj}. \How the second equation seem incorrect (and the claim that it follows from \cite[A. 29]{Gaj} seems a typo). To see this, note that subtracting the two equations of \cite[Cor. 12]{Gaj}
yields a very simple  equation in which $a,b$ are separated,
$\fr {S_a'(b)}{W_q'(b)}=d_a-a q$. The RHS   is a function of $a$ only (which equals $ c +m(a)-a q$ in our notation), and this seems different from $j(b)=k a$.

Such a separation did appear  in \eqr{J0a} in the exponential case (and yielded our   equation
\eqr{strEqb} which involves $b$ only), but  seems impossible  in the \me\ case, as  is apparent in  the RHS of \eqr{J0aPH}.
\eeR \fi 
%\input{MA}
\small
\bibliographystyle{amsalpha}
\bibliography{Pare37}

\end{document}